\pdfoutput=1
\documentclass{article}
\usepackage[vcentering]{geometry}
\usepackage{graphicx,xcolor}
\usepackage{tikz}
\usetikzlibrary{cd}
\usepackage[backend=biber,sorting=nyt]{biblatex}
\usepackage{amsmath,amssymb}
\usepackage{rsfso}
\allowdisplaybreaks[1]
\numberwithin{equation}{section}
\usepackage[bookmarksnumbered]{hyperref}
\usepackage{amsthm}
\theoremstyle{plain}
\newtheorem{theorem}{Theorem}[section]
\newtheorem{proposition}[theorem]{Proposition}
\newtheorem{lemma}[theorem]{Lemma}
\newtheorem{corollary}[theorem]{Corollary}
\theoremstyle{definition}
\newtheorem{definition}[theorem]{Definition}
\newtheorem{example}[theorem]{Example}
\newtheorem{remark}[theorem]{Remark}
\newcommand{\cadlag}{c{\`{a}}dl{\`{a}}g}

\newcommand{\caglad}{c{\`{a}}gl{\`{a}}d}

\newcommand{\lrbrack}{\mathopen{\rbrack}}
\newcommand{\rlbrack}{\mathclose{\lbrack}}
\newcommand{\Par}{\mathop{\mathrm{Par}}}
\usepackage{newtxtext}
\usepackage[varbb]{newtxmath}

\renewcommand{\thefootnote}{\fnsymbol{footnote}}

\begin{document}
\title{It\^{o}--F\"{o}llmer Calculus in Banach Spaces I: The It\^{o} Formula}
\author{Yuki Hirai}
\date{}
\maketitle

\begin{abstract}
    We prove F\"{o}llmer's pathwise It{\^o} formula
    for a Banach space-valued c\`{a}dl\`{a}g path.
    We also relax the assumption on the sequence of 
    partitions along which we treat the quadratic variation of a path.
\end{abstract}

\renewcommand{\thefootnote}{\arabic{footnote}}

\tableofcontents

\section{Introduction}
In his seminal paper~\cite{Foellmer_1981},
F\"{o}llmer presented a new perspective on
It{\^o}'s stochastic calculus.
The main theorem of F\"{o}llmer~\cite{Foellmer_1981}
states that a deterministic c\`{a}dl\`{a}g path satisfies the It{\^o} formula
provided it has quadratic variation
along a given sequence of partitions. 
This theorem enables us to construct the It{\^o} integral 
$\int_0^t f(X_{s-}) \mathrm{d}X_s$ for
a sufficiently nice function $f$ and a path $X$ with a quadratic variation.
This suggests the possibility 
of developing an analogue of the It{\^o} calculus
in completely analytic, probability-free situations.
We call this framework F\"{o}llmer's pathwise It{\^o} calculus
or, more simply, the It\^{o}--F\"{o}llmer calculus.
It can be regarded as a deterministic counterpart of the classical It{\^o} calculus.

Recently, the It\^{o}--F\"{o}llmer calculus has been receiving increasing attention 
from the viewpoint of its financial applications.
It is regarded as a useful tool to study financial theory under probability-free settings and has been used to construct financial strategies in a strictly pathwise manner 
(see, e.g., F\"{o}llmer and Schied~\cite{Foellmer_Schied_2013},
Schied~\cite{Schied_2014}, 
Davis, Ob{\l}{\'o}j and Raval~\cite{Davis_Obloj_Raval_2014},
and Schied, Speiser, and Voloshchenko~\cite{Schied_Speiser_Voloshchenko_2018}).
We expect that the It\^{o}--F\"{o}llmer calculus will
have a growing presence
in financial applications.

The It\^{o}--F\"{o}llmer calculus can be applied to
a stochastic process having quadratic variation.
A standard example of such a process is a semimartingale.
However, it is known that 
the class of processes possessing quadratic variation 
is strictly larger than that of semimartingales
(see, e.g., F\"{o}llmer~\cite{Foellmer_1981,Foellmer_1981b}).
In this sense, the It\^{o}--F\"{o}llmer calculus
enables us to extend stochastic integration theory beyond semimartingales.

There are several approaches to the pathwise construction of stochastic integration.
First, we mention classic studies by
Bichteler~\cite{Bichteler_1981},
Karandikar~\cite{Karandikar_1995},
and Willinger and Taqqu~\cite{Willinger_Taqqu_1988,Willinger_Taqqu_1989}, but
see also Nutz~\cite{Nutz_2012}
and {\L}ochowski~\cite{Lochowski_2014,Lochowski_2019}.
The theory of Vovk's outer measure and typical
paths was pioneered by Vovk~\cite{Vovk_2008,Vovk_2011,Vovk_2012,Vovk_2015,Vovk_2016}
and further developed by several authors, including 
Perkowski and Pr{\"o}mel~\cite{Perkowski_Proemel_2016,Perkowski_Proemel_2015b},
{\L}ochowski, Perkowski, and Pr{\"o}mel~\cite{Lochowski_Perkowski_Proemel_2018},
and Bartl, Kupper, and Neufeld~\cite{Bartl_Kupper_Neufeld_2019}.
Russo and Vallois~\cite{Russo_Vallois_1991,
    Russo_Vallois_1993,
    Russo_Vallois_1993b,
    Russo_Vallois_1995,
    Russo_Vallois_1996,
    Russo_Vallois_2000,
    Russo_Vallois_2007}
developed a theory called stochastic calculus via regularization.
The rough path theory,
pioneered by Lyons~\cite{Lyons_1998},
and its generalization have become 
important in stochastic calculus and its applications.
In addition, we refer to
Gubinelli~\cite{Gubinelli_2004},
Gubinelli and Tindel~\cite{Gubinelli_Tindel_2010},
Friz and Shekhar~\cite{Friz_Shekhar_2017}, 
and Friz and Zhang~\cite{Friz_Zhang_2018}.
Some studies have investigated the relation between the It\^{o}--F\"{o}llmer
calculus and rough path theory (see, e.g.,
Perkowski and Pr{\"o}mel~\cite{Perkowski_Proemel_2016},
and Friz and Hairer~\cite{Friz_Hairer_2014}).

Among the various pathwise methods,
we consider F\"{o}llmer's approach to be
the simplest and the most intuitively clear.
It needs only elementary arguments to establish 
calculation rules such as It{\^o}'s formula 
within this framework.
Moreover, the It\^{o}--F\"{o}llmer calculus requires only a minimal 
assumption that the integrator has quadratic variation.
We believe that these are advantages of F\"{o}llmer's theory,
and also that careful observation of this theory
helps us to understand the path properties of processes better
when we consider semimartingales and their stochastic integration.

Increasingly many works related to It\^{o}--F\"{o}llmer calculus
have been appearing recently.
First, we refer to Sondermann~\cite{Sondermann_2006},
Schied~\cite{Schied_2014}, Hirai~\cite{Hirai_2019}, and
Cont and Perkowski~\cite{Cont_Perkowski_2019}.
Schied~\cite{Schied_2016},
Mishura and Schied~\cite{Mishura_Schied_2016},
and Cont and Das~\cite{Cont_Das_2022}
construct deterministic continuous
paths with nontrivial quadratic variations.
See also Chiu and Cont~\cite{Chiu_Cont_2018}.
Functional extensions of the It\^{o}--F\"{o}llmer calculus 
have been developed by Dupire~\cite{Dupire_2009}, 
Cont and Fourni{\'e}~\cite{Cont_Fournie_2010,
Cont_Fournie_2010b,
Cont_Fournie_2013}, and 
Ananova and Cont~\cite{Ananova_Cont_2017}, for example.
Extension of the It\^{o}--F\"{o}llmer calculus
in terms of local times has been investigated in
Davis, Ob{\l}{\'o}j, and Raval~\cite{Davis_Obloj_Raval_2014},
Davis, Ob{\l}{\'o}j, and Siorpaes~\cite{Davis_Obloj_Siorpaes_2018},
{\L}ochowski et al.~\cite{Lochowski_Obloj_Promel_Siorpaes_2021b},
and Hirai~\cite{Hirai_2016,Hirai_2021}.

To our knowledge, however,
the It\^{o}--F\"{o}llmer calculus in an infinite dimensional setting
has not yet been sufficiently studied.
Stochastic integration in infinite dimensions naturally appears
when we treat stochastic partial differential equations
(see, e.g. Da Prato and Zabczyk~\cite{DaPrato_Zabczyk_2014}).
These have played an important role
in modelling term structures of interest rates or forward variances in mathematical finance,
and also in models of statistical mechanics and quantum field theories.
Then we aim to extend F\"{o}llmer's theory 
to Banach space-valued paths.
In this paper, we prove the It\^{o} formula for 
a path in a Banach space with a suitably defined quadratic variation.
We will study relations between various quadratic variations
and prove some transformation formulae for quadratic variations 
in our second paper in this series~\cite{Hirai_2022b}.
We not only generalize the state space of paths
but also relax the assumption on the sequence of
partitions along which we consider the quadratic variation.
In the context of the It\^{o}--F\"{o}llmer calculus,
two types of assumptions about a sequence of partitions are 
frequently used.
One is $\lvert \pi_n \rvert \to 0$, as used in F\"{o}llmer~\cite{Foellmer_1981},
and the other is the condition $O_t^{-}(X;\pi_n) \to 0$,
which is used in many papers handling continuous paths 
and some dealing with discontinuous paths
such as Vovk~\cite{Vovk_2015}.
In this paper, we introduce new conditions
to a sequence of partitions and a c\`{a}dl\`{a}g path
(Definition~\ref{2d}),
which gives a unified approach.

There have been many attempts to extend classical stochastic 
calculus to Banach or Hilbert space-valued processes.
Examples include Kunita~\cite{Kunita_1970},
Metivier~\cite{Metivier_1972},
Pellaumail~\cite{Pellaumail_1973},
Yor~\cite{Yor_1974},
Gravereaux and Pellaumail~\cite{Gravereaux_Pellaumail_1974},
Metivier and Pistone~\cite{Metivier_Pistone_1975},
Meyer~\cite{Meyer_1977a},
Metivier and Pellaumail~\cite{Metivier_Pellaumail_1980b},
Gy{\"o}ngy and Krylov~\cite{Gyongy_Krylov_1980,Gyongy_Krylov_1982},
Gy{\"o}ngy~\cite{Gyongy_1982},
Metivier~\cite{Metivier_1982},
Pratelli~\cite{Pratelli_1988},
Brooks and Dinculeanu~\cite{Brooks_Dinculeanu_1990},
Mikulevicius and Rozovskii~\cite{Mikulevicius_Rozovskii_1998,Mikulevicius_Rozovskii_1999},
Dinculeanu~\cite{Dinculeanu_2000},
De Donno and Pratelli~\cite{DeDonno_Pratelli_2006},
van Neerven, Veraar, and Weis~\cite{vanNeerven_Veraar_Weis_2007a},
Veraar and Yaroslavtsev~\cite{Veraar_Yaroslavtsev_2016},
and Yaroslavtsev~\cite{Yaroslavtsev_2020a}. 
Note that Di Girolami, Fabbri, and Russo~\cite{DiGirolami_Fabbri_Russo_2014}
treat quadratic covariation of Banach space-valued processes
within the framework of stochastic calculus via regularization,
with F\"{o}llmer's calculus in mind.

Our method can be interpreted as a deterministic counterpart
of these stochastic integration theories in Banach spaces.
Some of the works listed above,
such as Metivier and Pellaumail~\cite{Metivier_Pellaumail_1980b}
and Dinculeanu~\cite{Dinculeanu_2000},
give proofs of It\^{o}'s formula
in a similar manner to F\"{o}llmer's calculus.
One of the advantages of our approach 
appears in the statement of the It{\^o} formula.
For a function $f$ to satisfy the It{\^o} formula,
we require $f$ to be just $C^2$ class,
while a stochastic approach needs some additional assumptions about 
the boundedness of $f$ and its derivatives. 

Before explaining our contribution, 
we begin by summarizing the main result of F\"{o}llmer~\cite{Foellmer_1981}. 
Let $\Pi = (\pi_n)_{n \in \mathbb{N}}$ be a sequence of partitions of
$\mathbb{R}_{\geq 0}$ such that
$\lvert \pi_n \rvert \coloneqq \sup_{\lrbrack r,s \rbrack \in \pi_n} \lvert s-r \rvert$
tends to $0$ as $n \to \infty$.
We say that a {c\`{a}dl\`{a}g} path $X \colon \mathbb{R}_{\geq 0} \to \mathbb{R}$
has quadratic variation along $\Pi$ 
if there exists a {c\`{a}dl\`{a}g} increasing function $[X,X]$
such that for all $t \in \mathbb{R}_{\geq 0}$
\begin{enumerate}
    \item $\sum_{\lrbrack r,s \rbrack \in \pi_n} (X_{s \wedge t} - X_{r \wedge t})^2$
        converges to $[X,X]_t$ as $n \to \infty$, and
    \item $\Delta [X,X]_{t} = (\Delta X_t)^2$.
\end{enumerate}
An $\mathbb{R}^d$-valued {c\`{a}dl\`{a}g} path $X = (X^1,\dots,X^d)$
has quadratic variation along $\Pi$ if
the real-valued path $X^i + X^j$ has quadratic variation along the same sequence
for each $i$ and $j$.
F\"{o}llmer~\cite{Foellmer_1981} proved that if $X$ has quadratic variation, then for any $f \in C^2(\mathbb{R}^d)$
the path $t \mapsto f(X_t)$ satisfies It{\^o}'s formula.
That is,
\begin{align} \label{1.1c}
    f(X_t)
    & = 
    f(X_0)
    + \int_0^t \langle Df (X_{s-}), \mathrm{d}X_s \rangle
    + \frac{1}{2} \sum_{1 \leq i,j \leq d} 
        \int_0^t \frac{\partial^2 f}{\partial x_i \partial x_j}(X_{s-}) \mathrm{d}[X^i,X^j]^{\mathrm{c}}_s    \notag \\
    & \quad
    + \sum_{0 < s \leq t}
        \left\{ 
            \Delta f(X_s) 
            - \langle D f(X_{s-}),\Delta X_s \rangle
        \right\}
\end{align}
holds for all $t \in \mathbb{R}_{\geq 0}$.
The first term on the right-hand side of \eqref{1.1c}
is defined as the limit
\begin{equation*}
    \int_0^t \langle Df (X_{s-}), \mathrm{d}X_s \rangle
    =
    \lim_{n \to \infty} \sum_{\lrbrack r,s \rbrack \in \pi_n}
        \langle D f(X_r), X_{s \wedge t} - X_{r \wedge t} \rangle,
\end{equation*}
where $\langle \phantom{x},\phantom{x} \rangle$ denotes
the usual Euclidean inner product.
We call this limit the It\^{o}--F\"{o}llmer integral along $\Pi$.
F\"{o}llmer's theorem claims that if $X$ has quadratic variation along $\Pi$,
then the It\^{o}--F\"{o}llmer integral above
exists and it satisfies equation~\eqref{1.1c}.

As stated above, 
we aim to extend F\"{o}llmer's pathwise It{\^o} formula to Banach space-valued paths.
Let us describe a simplified version of our main result.
The precise statement will be given as Theorem~\ref{2e}
and Corollary~\ref{2fab}.
Let $E$ be a Banach space
and $E \widehat{\otimes}_{\alpha} E$ be the tensor product Banach space with respect to 
a reasonable crossnorm $\alpha$.
We say that an $E$-valued c\`{a}dl\`{a}g path $X$
has strong/weak $\alpha$-tensor quadratic variation along $\Pi = (\pi_n)$
if there is a c\`{a}dl\`{a}g path
${^\alpha [X,X]} \colon \mathbb{R}_{\geq 0} \to E \widehat{\otimes}_{\alpha} E$
of finite variation such that, for all $t \geq 0$,
\begin{enumerate}
    \item the sequence $\sum_{\lrbrack r,s \rbrack \in \pi_n} (X_{s \wedge t} - X_{r \wedge t})^{\otimes 2}$ 
        converges to ${^\alpha [X,X]}_t$ in the norm/weak topology of $E \widehat{\otimes}_{\alpha} E$, and
    \item the equation $\Delta {^\alpha [X,Y]}_t = \Delta X_t^{\otimes 2}$ holds.
\end{enumerate}
The path $X$ has finite 2-variation along $\Pi$ if
\begin{equation*}
    V^2(X;\Pi)_t
    \coloneqq 
    \sup_{n \in \mathbb{N}} 
        \sum_{\lrbrack r,s \rbrack \in \pi_n} \lVert X_{s \wedge t} - X_{r \wedge t} \rVert^2
    < \infty
\end{equation*}
for all $t \geq 0$. 
We say that a sequence of partitions $(\pi_n)$
satisfies Condition (C) for a path $X \colon \lbrack 0,\infty \rlbrack \to E$
in a Banach space if it satisfies Conditions~(C1)--(C3) of Definition~\ref{2d}.
Roughly speaking, Conditions~(C1) and (C2) state that 
$(\pi_n)$ reconstructs the information of the jumps of $X$.
Condition~(C3) means that 
$(\pi_n)$ controls the oscillation of $X$ in some sense.
Under these settings, our main result, the It\^{o} formula, is stated as follows.
Let $X \colon \mathbb{R}_{\geq 0} \to E$
be a c\`{a}dl\`{a}g path that has strong/weak $\alpha$-tensor quadratic variation 
and finite 2-variation along $(\pi_n)_{n \in \mathbb{N}}$,
and let $A \colon \mathbb{R}_{\geq 0} \to F$ be a c\`{a}dl\`{a}g path
of finite variation in a Banach space.
Suppose that $(\pi_n)$ satisfies Condition (C) for $(A,X)$ and 
the left-side discretization of $(A,X)$ along $(\pi_n)$
approximates $(A_{-},X_{-})$ pointwise
(see Definition~\ref{2b} for the exact definition).
If $f \colon F \times E \to G$ is a $C^{1,2}$ function such that the second derivative
induces a continuous map $D^2_x f \colon F \times E \to \mathcal{L}(E \widehat{\otimes}_{\alpha} E,G)$,
then the composite function $f(A,X)$ satisfies
\begin{multline*}
    f(A_t,X_t) - f(A_0,X_0)
    = 
    \int_{0}^{t} D_a f(A_{s-},X_{s-}) \mathrm{d}A^{\mathrm{c}}_s  
    + \int_0^t D_x f (A_{s-},X_{s-}) \mathrm{d}X_s      \\ 
    + \frac{1}{2} \int_{0}^{t} D_x^2 f(A_{s-},X_{s-})\mathrm{d}{^\alpha [X,X]}^{\mathrm{c}}_s  
    + \sum_{0 < s \leq t} \left\{ \Delta f(A_{s},X_{s}) - D_x f(A_{s-},X_{s-})\Delta X_s \right\}. 
\end{multline*}
The second integral on the right-hand side is
defined respectively as the strong/weak limit of left-side Riemannian sums along $(\pi_n)$.

To conclude this section,
we give an outline of the remainder of this paper.
Section~2 is a preliminary part of this article.
We introduce basic notation and terminology in the first subsection.
The next subsection is devoted to a review of 
c\`{a}dl\`{a}g paths and Stieltjes integrals in Banach spaces.
In Section~3, we set up basic notions in It\^{o}--F\"{o}llmer calculus
in Banach spaces and state the main results of the paper 
(Theorem~\ref{2e} and Corollary~\ref{2fab}).
In Section~4, we study conditions
on the sequence of partitions and the relation between them and 
c\`{a}dl\`{a}g paths.
Fundamental properties of quadratic variations are studied in Section~5.
The purpose of Section~6 is to show Lemma~\ref{2g}, 
which is essentially used in the proof of the main theorem.
In the last section of the main part, Section~7, 
we prove the It{\^o} formula for a Banach space-valued path having quadratic variation.
In appendices, we present some auxiliary results
related to differential calculus and integration in Banach spaces.

\section{Preliminaries}

\subsection{Notations and terminologies}

In this section, we introduce basic notation and terminology used throughout this article.

The symbol $\mathbb{N}$ denotes the set of natural numbers $\{ 0,1,2,\dots \}$
and $\mathbb{R}$ denotes that of real numbers.
If $A$ is a subset of $\mathbb{R}$ and $a \in \mathbb{R}$, we define 
$A_{\geq a} = \{ x \in A \mid x \geq a \}$.

If $E$ and $F$ are two real Banach spaces,
$\mathcal{L}(E,F)$ denotes the space of bounded linear maps from $E$ to $F$.
In addition, given another Banach space $G$,
we define $\mathcal{L}^{(2)}(E,F;G)$ as the 
space of bounded bilinear maps from $E \times F$ to $G$.
Recall that $\mathcal{L}(E,F)$ and $\mathcal{L}^{(2)}(E,F;G)$ are Banach spaces with norms
\begin{equation*}
    %1
    \lVert L \rVert 
    = \sup_{\lVert x \rVert \neq 0} \frac{\lVert Lx \rVert}{\lVert x \rVert}, \qquad
    %2
    \lVert B \rVert 
    = \sup_{\lVert x \rVert, \lVert y \rVert \neq 0} 
        \frac{\lVert B(x,y) \rVert}{\lVert x \rVert \lVert y \rVert},
\end{equation*}
respectively.

Now we introduce another topology on the space $\mathcal{L}(E,G)$.
Let $\mathcal{K}_E$ be the family of all compact subsets of $E$.
For each $K \in \mathcal{K}_E$, define a seminorm $\rho_K$
by the formula 
\begin{equation} \label{2.1b}
    \rho_K(L) = \inf \{ C > 0 \mid \forall x \in K, \ \lVert Lx \rVert_F \leq C \lVert x \rVert_E \}
\end{equation}
for each $L \in \mathcal{L}(E,F)$.
Then the family $(\rho_K)_{K \in \mathcal{K}_E}$ induces a locally convex Hausdorff 
topology on $\mathcal{L}(E,F)$.
We use the symbol $\mathcal{L}_{\mathrm{c}}(E,F)$ for this topological vector space.

Let $[0,\infty \rlbrack = \mathbb{R}_{\geq 0} = \{ r \in \mathbb{R} \mid r \geq 0 \}$
and let $E$ be a Hausdorff topological vector space.
A \emph{c\`{a}dl\`{a}g} path in $E$
is a function $X \colon \mathbb{R}_{\geq 0} \to E$
that is right continuous at every $t \geq 0$
and has a left limit at every $t > 0$.
The terms \emph{RCLL} and \emph{right-regular} are also used
to stand for the same property.
Similarly, a \emph{c\`{a}gl\`{a}d} (also called \emph{LCRL} or \emph{left-regular}) path in $E$
is a function $X \colon \mathbb{R}_{\geq 0} \to E$
that is left continuous on $\lrbrack 0,\infty \rlbrack$
and has right limits on $[0,\infty \rlbrack$. 
The symbols $D([0,\infty \rlbrack, E)$
and $D(\mathbb{R}_{\geq 0},E)$ denote the set of 
all c\`{a}dl\`{a}g paths in $E$.
If $X$ is an element of $D(\mathbb{R}_{\geq 0}, E)$, we define
\begin{equation*}
    X(t-) = \lim_{s \uparrow\uparrow t} X(s) = \lim_{s \to t, s < t} X_s, \qquad
    \Delta X(t) = X(t) - X(t-).
\end{equation*}
We also use $X_t$, $X_{t-}$, and $\Delta X_t$ 
to indicate the values $X(t)$, $X(t-)$, and $\Delta X(t)$, respectively.
Next, set
\begin{gather*}
    %1
    D(X) 
        = \{ t \in \mathbb{R}_{\geq 0} \mid \lVert \Delta X_t \rVert \neq 0 \},  \\
    %2
    D_{\varepsilon}(X) 
        = \{ t \in \mathbb{R}_{\geq 0} \mid \lVert \Delta X_t \rVert \geq \varepsilon \},  \\
    %3
    D^{\varepsilon}(X) = D(X) \setminus D_{\varepsilon}(X)
        = \{ t \in \mathbb{R}_{\geq 0} \mid 0 < \lVert \Delta X_t \rVert < \varepsilon \}.
\end{gather*}
We simply write $D$, $D_{\varepsilon}$, and $D^{\varepsilon}$ if there is no ambiguity.
Given a discrete set $D \subset [0,\infty[$ and a c\`{a}dl\`{a}g path $X$, we define
\begin{equation*}
    J_D(X)_t = J(D;X)_t = \sum_{0 < s \leq t} \Delta X_s 1_{D}(s).
\end{equation*}
Recall that a subset $D$ is discrete if it is a discrete topological subspace of 
$\lbrack 0,\infty \rlbrack$; i.e. every element of $D$ is
an isolated point with respect to the subspace topology.
By assumption, the set $D \cap [0,t]$ is finite for all $t$, and therefore
the summation above is well-defined.
Then $J_D(X)$ is a c\`{a}dl\`{a}g path of finite variation.
For abbreviation, we often write
$J_\varepsilon(X)$ instead of $J(D_{\varepsilon}(X);X)$.

Throughout this paper,
the term \emph{partition of $\mathbb{R}_{\geq 0}$} always 
means the set of intervals of the form
$\pi = \{ \lrbrack t_i,t_{i+1} \rbrack ; i \in \mathbb{N} \}$
which satisfies $0=t_0 < t_1 < \dots \to \infty$. 
The set of all partitions of $\mathbb{R}_{\geq 0}$ is 
denoted by $\Par(\mathbb{R}_{\geq 0})$ or $\Par (\lbrack 0,\infty \rlbrack )$.
Similarly, given a compact interval $[a,b] \subset \mathbb{R}_{\geq 0}$,
let $\Par([a,b])$ be the set of all finite partitions of the 
form $\pi = \{ \lrbrack t_i,t_{i+1} \rbrack ; 0 \leq i \leq n-1 \}$
with $a=t_0 < t_1 < \dots < t_n = b$.
For a partition $\pi$ of $\mathbb{R}_{\geq 0}$ or a compact interval,
we define $\pi^{\mathrm{p}} \subset \mathbb{R}_{\geq 0}$
as the set of all endpoints of elements of $\pi$.
If $\pi = \{ \lrbrack t_i,t_{i+1} \rbrack ; i \in I \}$,
then $\pi^{\mathrm{p}} = \{ t_i, t_{i+1} ; i \in I \}$. 

\subsection{Remarks on c\`{a}dl\`{a}g paths and Stieltjes integration}

In this subsection, we review some basic properties of c\`{a}dl\`{a}g paths that will be 
referred to later.

Let $E$ be a Banach space.
$E$-valued right continuous and left continuous step functions are functions of the form
\begin{equation*}
    \sum_{i \in \mathbb{N}} 1_{[t_i,t_{i+1}\mathclose{[}} a_{i}, \qquad
    1_{\{ 0 \}} b_0  + \sum_{i \in \mathbb{N}} 1_{\lrbrack t_i,t_{i+1} \rbrack} b_{i+1}, 
\end{equation*}
respectively, where $0 = t_0 < t_1 < \dots < t_n < \dots \to \infty$
and $a_i$, $b_i \in E$ for all $i \in \mathbb{N}$.
Right continuous step functions are c\`{a}dl\`{a}g
and left continuous step functions are c\`{a}gl\`{a}d.
Every right continuous step function $f= \sum_{i \in \mathbb{N}} 1_{[t_i,t_{i+1}\mathclose{[}} a_{i}$
is strongly $\mathcal{B}(\mathbb{R}_{\geq 0})/\mathcal{B}(E)$ measurable
because it is the pointwise limit of the sequence $(f_n)$ defined by
$f_n = \sum_{0 \leq i \leq n} 1_{[t_i,t_{i+1}\mathclose{[}} a_{i}$.

A c\`{a}dl\`{a}g path in a Banach space satisfies the following properties.

\begin{lemma} \label{3b}
Let $f$ be a c\`{a}dl\`{a}g path in a Banach space $E$.
\begin{enumerate}
    \item For every $C > 0$, there are only finitely many $s$
        satisfying $\lVert \Delta f_{s} \rVert_E > C$
        in each compact interval of $[0,\infty \rlbrack$.
    \item The image $f(I)$ of any compact interval $I \subset [0,\infty \rlbrack$ is relatively compact in $E$.
    \item Suppose that every jump of $f$ is smaller than $C > 0$
        on a compact interval $I \subset [0,\infty \rlbrack$.
        Then for all $\varepsilon > 0$, there exists a $\delta > 0$
        such that $\lVert f(s) - f(u) \rVert_E < C + \varepsilon$ holds
        for any $s,u \in I$ satisfying $\lvert s - u \rvert < \delta$.
    \item The path $f$ is the uniform limit of some sequence
        of right continuous step functions on each compact interval.
    \item For every $t > 0$ and $\varepsilon > 0$,
        there is a partition $\pi \in \Par([0,t])$ that satisfies
        \begin{equation*}
            \omega(f, \lbrack r,s \rlbrack) \coloneqq \sup_{u,v \in \lbrack r,s \rlbrack} \lVert f(u) - f(v) \rVert < \varepsilon
        \end{equation*}
        for all $\lrbrack r,s \rbrack \in \pi$.
\end{enumerate}
\end{lemma}

One can find an analogue of Lemma~\ref{3b}
for c\`{a}dl\`{a}g paths in an arbitrary separable complete metric space
in Billingsley~\cite[122]{Billingsley_1999}.

Next, recall that a function $f \colon [0,\infty \rlbrack \to E$ is
of bounded variation on a compact interval $I \subset [0,\infty \rlbrack$ if
\begin{equation*}
    V(f;I)
    \coloneqq
    \sup_{\pi \in \Par I}
        \sum_{\lrbrack r,s \rbrack \in \pi} \lVert f(s) - f(r) \rVert_E < \infty.
\end{equation*}
For convenience, set $V(f;\emptyset) = 0$ and $V(f;[a,a]) = 0$
for every $a \in [0,\infty \rlbrack$.
The function $f$ has finite variation 
if it has bounded variation on every compact subinterval of $[0,\infty \rlbrack$.
The set of all c\`{a}dl\`{a}g paths of finite variation 
in $E$ is denoted by $FV(\mathbb{R}_{\geq 0},E)$
or $FV(\lbrack 0,\infty \rlbrack,E)$.
We define the total variation path $V(f)$ of a function $f$
of $FV(\mathbb{R}_{\geq 0},E)$ by $V(f)_t = V(f;[0,t])$.
Then $V(f)$ is increasing and satisfies $V(f)_0 = 0$.

We list several basic properties of a path of finite variation below.
See Dinculeanu~\cite[{\S}18]{Dinculeanu_2000} for a proof.

\begin{lemma} \label{3c}
Let $f \colon [0,\infty \rlbrack \to E$ be a c\`{a}dl\`{a}g path of finite variation in a Banach space.
\begin{enumerate}
    \item If $a,b,c \in [0,\infty \rlbrack$ satisfy $a \leq b \leq c$, then
        \begin{equation*}
            V(f,[a,c]) = V(f,[a,b]) + V(f,[b,c]).
        \end{equation*}
    \item The total variation path $V(f)$ is c\`{a}dl\`{a}g.
    \item The jump of $V(f)$ at $t \geq 0$ is given by $\Delta V(f)(t) = \lVert \Delta f(t) \rVert_E$.
    \item The family $(\Delta f(s))_{s \in [0,t]}$ is absolutely summable for all $t \geq 0$; i.e.
        \begin{equation*}
            \sup \left\{ 
                \sum_{s \in F} \lVert \Delta f(s) \rVert_E \, 
                \middle\vert \,
                \text{$F$: a finite subset of $[0,t]$} \right\} < \infty.
        \end{equation*}
    \item The function $f^{\mathrm{d}}$ defined by
        \begin{equation*}
            f^{\mathrm{d}}(t) = \sum_{0 < s \leq t} \Delta f(s).
        \end{equation*}
        is again a c\`{a}dl\`{a}g path of finite variation.
\end{enumerate}
\end{lemma}

Note that the summation in (v) of Lemma~\ref{3c}
is defined in the following manner.
Let $D$ be the set of all finite subsets of $\lrbrack 0,t \rbrack$.
We regard $D$ as a directed set with the order defined by inclusion.
Then the net $(\sum_{s \in d} \Delta f(s))_{d \in D}$ 
converges in $E$ by Condition (iv) of Lemma~\ref{3c}, 
and hence we can define
\begin{equation*}
    \sum_{0 < s \leq t} \Delta f(s) = \lim_{d} \sum_{s \in d} \Delta f(s).
\end{equation*}
The function $f^{\mathrm{d}}$ defined in Proposition~\ref{3c}
is called the discontinuous part of $f$.
We also define $f^{\mathrm{c}} = f - f^{\mathrm{d}}$
and call it the continuous part of $f$.

Let $\mathcal{I}$ the semiring of subsets of $\mathbb{R}_{\geq 0}$
consisting of all bounded intervals of the form $\lrbrack a,b \rbrack$
and the singleton $\{ 0 \}$.
Given an $f \in D(\mathbb{R}_{\geq 0},E)$, define
\begin{equation*}
    \mu_{f}(\lrbrack a,b \rbrack) = f(b) - f(a), \quad 
    \mu_{f}(\{ 0 \}) = f(0)
\end{equation*}
for any two real numbers satisfying $0 \leq a \leq b$.
If $f$ has finite variation,
the function $\mu_f \colon \mathcal{I} \to E$
can be uniquely extended to a $\sigma$-additive 
measure defined on the $\delta$-ring $\mathcal{D}$
generated by $\mathcal{I}$.
Refer to Theorem~18.19 of Dinculeanu~\cite[208]{Dinculeanu_2000}
for a proof.
Notice that this correspondence between a function and a measure
satisfies $(\mu_f)^{\mathrm{d}} = \mu_{f^\mathrm{d}}$ and 
$(\mu_f)^{\mathrm{c}} = \mu_{f^{\mathrm{c}}}$.
See Appendix B.1 for the definition of
the measures $(\mu_f)^{\mathrm{d}}$ and $(\mu_f)^{\mathrm{c}}$.

Because there is a measure $\mu_f$ associated with $f$,
we can consider the Stieltjes integral with respect to $f$.
Let $B \colon F \times E \to G$ be a continuous bilinear map
between Banach spaces.
Set $L^{1}_{\mathrm{loc}}(\mu_f;F) = L^{1}_{\mathrm{loc}}(\lvert \mu_f \rvert;F)$,
where $\lvert \mu_f \rvert$ denotes the variation measure of $\mu_f$
introduced in Appendix B.1.
Then for each $g \in L^{1}_{\mathrm{loc}}(\mu_f;F)$ and
a bounded interval $I = \lrbrack a,b \rbrack$,
the Stieltjes integral is defined by the formula
\begin{equation*}
    \int_a^b B(g(s),\mathrm{d}f(s))
    =
    \int_{I} B(g(s), \mathrm{d}f(s))
    \coloneqq
    \int_{I} B(g(s), \mu_f(\mathrm{d}s)).
\end{equation*}
See Appendix B.1 for the definition of $\int_{I} B(g(s), \mu_f(\mathrm{d}s))$.
If $g \colon \mathbb{R}_{\geq 0} \to \mathcal{L}_{\mathrm{c}}(E,G)$ is a {\caglad} path,
we can also define the Stieltjes integral as
\begin{equation*}
    \int_{I} g(s)\, \mathrm{d}f(s)
    =
    \int_{I} g(s)\, \mu_f(\mathrm{d}s),
\end{equation*}
where the integral on the right-hand side is constructed in Appendix B.2. 
Finally, note that the decomposition $f = f^{\mathrm{c}} + f^{\mathrm{d}}$
gives a decomposition of the integral
\begin{equation*}
    \int_{I} B(g(s), \mathrm{d}f(s))
    =
    \int_{I} B(g(s), \mathrm{d}f^{\mathrm{c}}(s))
    + \int_{I} B(g(s), \mathrm{d}f^{\mathrm{d}}(s)).
\end{equation*}

\section{Settings and the main result}

In this section, 
we introduce the main results of this paper,
namely, the It{\^o} formula
within the framework of the It\^{o}--F\"{o}llmer calculus in Banach spaces.
The statement is to be found in Theorem~\ref{2e} and Corollary~\ref{2fab}.

First, we introduce some notations about difference operators.
Given a function $X \colon \mathbb{R}_{\geq 0} \to E$
and $t \geq 0$,
define functions $\delta X$ and $\delta X_t$ on $\mathcal{I}$ into $E$ by the formulae
\begin{equation*}
    \delta X(I) = \delta X(\lrbrack r,s \rbrack) = X_{s} - X_{r}, \qquad
    \delta X_t(I) = \delta X(I \cap [0,t]) = X_{s \wedge t} - X_{r \wedge t}
\end{equation*}
for each $I = \lrbrack r,s \rbrack \in \mathcal{I}$.
We also define $\delta X(\{ 0 \}) = \delta X_t(\{ 0 \}) = X_0$.
Moreover, consider a bilinear map $B \colon F \times E \to G$ in Banach spaces
and another function $Y \colon \mathbb{R}_{\geq 0} \to F$.
Then we define functions $B(Y,\delta X)$ and $B(Y,\delta X_t)$
from $\mathcal{I}$ to $G$ by the formulae
\begin{equation*}
    B(Y,\delta X)(\lrbrack r,s \rbrack) = B(Y_r,\delta X(\lrbrack r,s \rbrack)), \qquad 
    B(Y, \delta X_t)(\lrbrack r,s \rbrack) = B(Y_r,\delta X(\lrbrack r,s \rbrack \cap [0,t]))
\end{equation*}
for $I = \lrbrack r,s \rbrack \in \mathcal{I}$
and
\begin{equation*}
    B(Y,\delta X)(\{ 0 \}) = B(Y,\delta X_t)(\{ 0 \}) = B(Y_0,X_0).
\end{equation*}
By using this notation, the left-side Riemannian sum along a partition $\pi$ is expressed as
\begin{equation*}
    \sum_{\lrbrack r,s \rbrack \in \pi} B(Y_r, X_{t \wedge s} - X_{t \wedge r})
    =
    \sum_{\lrbrack r,s \rbrack \in \pi} B(Y_{t \wedge r}, X_{t \wedge s} - X_{t \wedge r})
    =
    \sum_{I \in \pi} B(Y,\delta X_t)(I).
\end{equation*}
We also define $B(\delta Y,X)$ and $B(\delta Y_t,X)$ in a similar manner.
The function $B(\delta X_t,\delta Y_t) \colon \mathcal{I} \to G$
is defined as the composition of $B$ and $(\delta X_t, \delta Y_t) \colon \mathcal{I} \to F \times E$.

\begin{definition} \label{2b}
Let $E$, $F$, and $G$ be Banach spaces and 
let $B \colon E \times F \to G$ be a bounded bilinear map.
Suppose that $\Pi = (\pi_n)_{n \in \mathbb{N}}$ is a sequence of partitions of $\mathbb{R}_{\geq 0}$.
\begin{enumerate}
    \item A path $(X,Y) \in D(\mathbb{R}_{\geq 0},E \times F)$ has
        \emph{$B$-quadratic covariation along $\Pi$}
        if there exists a $G$-valued c\`{a}dl\`{a}g path $Q_B(X,Y)$ of 
        finite variation satisfying the following conditions:
        \begin{enumerate}
            \item for all $t \in \mathbb{R}_{\geq 0}$
                \begin{equation*}
                    \lim_{n \to \infty} \sum_{I \in \pi_n} B(\delta X_t, \delta Y_t)(I) 
                    =
                    Q_B(X,Y)_t
                \end{equation*}
                in the norm topology of $G$;
            \item for all $t \in \mathbb{R}_{\geq 0}$, the jump of $Q_B(X,Y)$ is given by
                \begin{equation*}
                    \Delta Q_B(X,Y)_t = B(\Delta X_t,\Delta Y_t).
                \end{equation*}
        \end{enumerate}
        Then the path $Q_B(X,Y)$ is called
        the $B$-quadratic covariation of $X$ and $Y$.
        If $E=F$ and $X=Y$,
    we call $Q_B(X,X)$ the $B$-quadratic variation of $X$.
    \item If the convergence of (i)-(a) holds in the weak topology of $G$,
        we say that $(X,Y)$ has the
        \emph{weak $B$-quadratic covariation} $Q_B(X,Y)$.
\end{enumerate}
\end{definition}

We often call a $B$-quadratic covariation a \emph{strong} $B$-quadratic covariation
if we stress that the convergence holds in the norm topology.
The quadratic covariation $Q_B(X,Y)$ depends
on the sequence of partitions $\Pi$.
Given a partition $\pi$, we often write
\begin{equation*}
    Q_B^{\pi}(X,Y)_t = \sum_{I \in \pi} B(\delta X_t,\delta Y_t)(I).
\end{equation*}
Because $X$ and $Y$ are c\`{a}dl\`{a}g, 
the map $t \mapsto Q^{\pi}_B (X,Y)_t$
is also c\`{a}dl\`{a}g. It is not, however, of finite variation
unless $X$ and $Y$ are of finite variation.
With this notation, we can say that the
strong/weak $B$-quadratic covariation $Q_B(X,Y)$ is 
the pointwise limit of $(Q_B^{\pi_n}(X,Y))_{n \in \mathbb{N}}$
respectively in the norm/weak topology
satisfying the jump condition (i)-(b) of Definition~\ref{2b}.

An important class of a bounded bilinear map $B$ is the 
canonical bilinear map into a tensor product of Banach spaces.
Let $\alpha$ be a norm on the algebraic tensor product $E \otimes F$
of two Banach spaces.
The norm $\alpha$ is called a reasonable crossnorm if it satisfies 
the following two conditions:
\begin{enumerate}
    \item the inequality $\alpha(x \otimes y) \leq \lVert x \rVert \lVert y \rVert$
        holds for all $x \in E$ and $y \in F$;
    \item the inequality $\lVert x^* \otimes y^* \rVert_{(E \otimes F,\alpha)^*} \leq \lVert x^* \rVert \lVert y^* \rVert$
        holds for all $x^* \in E^*$ and $y^* \in F^*$,
        where $\lVert \phantom{x} \rVert_{(E \otimes F,\alpha)^*}$ denotes the usual operator norm
        on the normed space $(E \otimes F,\alpha)$.
\end{enumerate}
The completion of the normed space $(E \otimes F,\alpha)$,
which is generally incomplete, is denoted by $E \widehat{\otimes}_{\alpha} F$.
See Diestel and Uhl~\cite{Diestel_Uhl_1977} and Ryan~\cite{Ryan_2002}
for basic facts about tensor products of Banach spaces.
The quadratic covariation of $(X,Y)$ with respect to the canonical bilinear map
$\otimes \colon E \times F \to E \widehat{\otimes}_{\alpha} F$
is denoted by ${^\alpha [X,Y]}$,
and it is called the \emph{$\alpha$-tensor quadratic covariation}.
We also write $[X,Y]^{\pi} = Q_{\otimes}^{\pi}(X,Y)$
and call it the discrete tensor quadratic covariation of $(X,Y)$ along $\pi$.
If $E = F$, we can consider $\alpha$-tensor quadratic covariations
${^\alpha [X,Y]}$ and ${^\alpha [Y,X]}$. 
Although one of them exists if and only if the other does,
they are not equal in general.
The path ${^\alpha [X,X]}$ is called the \emph{$\alpha$-tensor quadratic variation}.

There are various important reasonable crossnorms in Banach space theory.
The greatest crossnorm $\gamma$, also called the projective tensor norm
\footnote{
    The projective tensor norm is often denoted by $\pi$.
    We use $\gamma$, following Diestel and Uhl~\cite{Diestel_Uhl_1977},
    because the symbol $\pi$ is used to indicate a partition in this article.
},
is defined by the following formula
\begin{equation*}
    \gamma(z) 
    = 
    \inf \left\{ \sum_{i \in I} \lVert x_i \rVert_E \lVert y_i \rVert_F \,
        \middle\vert \,
        \text{$I$ is finite, 
            $x_i \in E$ and $y_i \in F$ for all $i \in I$,
        and 
        $\displaystyle z = \sum_{i \in I} x_i \otimes y_i$}
    \right\}.
\end{equation*}
We simply write $E \widehat{\otimes} F = E \widehat{\otimes}_{\gamma} F$
and call it the projective tensor product of $E$ and $F$.
$\gamma$-tensor quadratic variations are also called \emph{projective tensor quadratic variations}.
An advantage of the projective tensor product is that there is 
a canonical isometric isomorphism
$
    \mathcal{L}(E,\mathcal{L}(E,G)) 
    \simeq \mathcal{L}(E \widehat{\otimes} E,G) 
    \simeq \mathcal{L}^{(2)}(E,E;G)
$.
We use this identification without mention.
If $E$ and $F$ are Hilbert spaces, Hilbert--Schmidt crossnorm is also 
a natural object.
For $L^p(\mu_1)$ and $L^p(\mu_2)$ defined on some measure spaces, 
there is a crossnorm $\triangle_p$ such that
$L^p(\mu_1) \widehat{\otimes}_{\triangle_p} L^p(\mu_2) \cong L^p(\mu_1 \otimes \mu_2)$.
See Defant and Floret~\cite[Section 7]{Defant_Floret_1993}.

The tensor quadratic variation of an $\mathbb{R}^d$-valued path
$X= (X^1,\dots,X^d)$ has the matrix representation
\begin{equation*}
    [X,X]_t
    =
    \begin{pmatrix}
        [X^1,X^1]_t & \cdots & [X^1,X^d]_t    \\
        \vdots & \ddots & \vdots   \\
        [X^d,X^1]_t & \cdots & [X^d,X^d]_t
    \end{pmatrix}.
\end{equation*}
A c\`{a}dl\`{a}g path $X \colon \mathbb{R}_{\geq 0} \to \mathbb{R}^d$
has tensor quadratic variation if and only if
it has quadratic variation in the sense of 
Definition~2.3 of Hirai~\cite{Hirai_2019}.

\begin{remark} \label{2bb}
\begin{enumerate}
    \item Our definition of strong tensor quadratic variations can be regarded
        as a pathwise analogue of tensor quadratic variations in classical stochastic integration theory
        in infinite dimensions such as Metivier and Pellaumail~\cite{Metivier_Pellaumail_1980b}
        and Metivier~\cite{Metivier_1982}.
        See also Dinculeanu~\cite{Dinculeanu_2000}.
        Although one can define tensor quadratic variations in general Banach spaces,
        classical existence results only deal with Hilbert spaces.
    \item Another important approach is developing recently in the context of
        the martingale theory in UMD Banach spaces.
        Yaroslavtsev~\cite{Yaroslavtsev_2020a} shows that a local martingale in a UMD Banach space
        has the covariation bilinear form $\lbrack\!\lbrack M \rbrack\!\rbrack$.
        In our terminology, the covariation bilinear form corresponds to
        the cylindrical quadratic variation in the sense of
        Corollary 3.5 in Hirai~\cite{Hirai_2022b}.
        By this corollary, we see that a {\cadlag} path with weak tensor quadratic variation 
        has cylindrical quadratic variation. 
\end{enumerate}
\end{remark}

Now we introduce a different type of quadratic variation,
namely, scalar quadratic variation.
Again we assume that $\Pi = (\pi_n)$ is a sequence of partitions of
$\mathbb{R}_{\geq 0}$.

\begin{definition} \label{2c}
Let $E$ be a Banach space and $X$ be an $E$-valued c\`{a}dl\`{a}g path.
\begin{enumerate}
    \item The \emph{$2$-variation of $X$ on $[0,t]$ along $\Pi$} is defined by
        \begin{equation*}
            V^2(X;\Pi)_t = \sup_{n \in \mathbb{N}} \sum_{I \in \pi_n} \lVert \delta X_t(I) \rVert^2.
        \end{equation*}
        We say the $X$ has finite $2$-variation along $\Pi$ if $V^2(X;\Pi)_t < \infty$
        for all $t \geq 0$.
    \item A c\`{a}dl\`{a}g path $X \colon \mathbb{R}_{\geq 0} \to E$ has
        \emph{scalar quadratic variation along $\Pi$}
        if there exists a real-valued c\`{a}dl\`{a}g increasing path
        $Q(X)$ such that
        \begin{enumerate}
            \item for all $t \in \mathbb{R}_{\geq 0}$,
                \begin{equation*}
                    \sum_{I \in \pi_n} \lVert \delta X_t(I) \rVert^2 
                    \xrightarrow[n \to \infty]{} Q(X)_t,
                \end{equation*}
            \item for all $t \in \mathbb{R}_{\geq 0}$, the jump of $Q(X)$ at $t$ is given by
                $\Delta Q(X)_t = \lVert \Delta X_t \rVert_E^2$.
        \end{enumerate}
        We call the increasing path $Q(X)$ the scalar quadratic variation of $X$
        along $(\pi_n)$.
\end{enumerate}
\end{definition}

If $X$ has scalar quadratic variation along $\Pi$,
then it has finite $2$-variation along $\Pi$.

The scalar quadratic variation of a Hilbert space valued path 
$X$ coincides with the bilinear quadratic variation $Q_{\langle\phantom{x},\phantom{x}\rangle}(X,X)$,
where $\langle\phantom{x},\phantom{x}\rangle$ is the inner product of the state space.
If a c\`{a}dl\`{a}g path $X = (X^1,\dots,X^d) \colon \mathbb{R}_{\geq 0} \to \mathbb{R}^d$
has tensor quadratic variation along $(\pi_n)$,
then it has scalar quadratic variation given by
\begin{equation*}
    Q(X)_t = \mathop{\mathrm{Trace}} [X,X]_t = \sum_{1 \leq i \leq d} [X^i,X^i]_t.
\end{equation*}
This trace representation is still valid for 
Hilbert space-valued c\`{a}dl\`{a}g paths.
This result is proved in Hirai~\cite{Hirai_2022b}.

Next, we introduce some conditions
on a sequence of partitions and a c\`{a}dl\`{a}g path.
Let $\pi \in \Par(\lbrack 0,\infty \rlbrack)$ and $t \in \lrbrack 0,\infty \mathclose{[}$.
The symbol $\pi(t)$ denotes the element of $\pi$ that contains $t$.
By definition, there exists only one such interval.
In addition, set
$\overline{\pi}(t) = \sup \pi(t)$ and $\underline{\pi}(t) = \inf \pi(t)$.
Then we have $\pi(t) = \lrbrack \underline{\pi}(t),\overline{\pi}(t) \rbrack$.

Let $f \colon S \to E$ be a function into a Banach space and $A$ be a subset of $S$.
Define the oscillation of $f$ on $A$ by
\begin{equation*}
    \omega(f;A) = \sup_{x,y \in A} \lVert f(x) - f(y) \rVert_E.
\end{equation*}
Using this notation, we introduce two kinds of oscillation of
a path $X \in D(\mathbb{R}_{\geq 0},E)$ 
along a partition $\pi \in \Par(\mathbb{R}_{\geq 0})$ as follows.
\begin{gather*}
    O^+_{t}(X;\pi) = \sup_{\mathopen{]}r,s] \in \pi} \omega(X;\lrbrack r,s \rbrack \cap [0,t]),  \\
    O^{-}_t(X;\pi)
    =
    \sup_{\mathopen{]}r,s] \in \pi} \omega(X;\mathopen{]}r,s\mathclose{[} \cap [0,t])
    =
    \sup_{\mathopen{]}r,s] \in \pi} \omega(X;\mathopen{[}r,s\mathclose{[} \cap [0,t]).
\end{gather*}
The second equality in the definition of $O^{-}_t(X;\pi)$
is valid because $X$ is supposed to be right continuous.
These oscillations satisfy the relation $O_t^-(X;\pi) \leq O_t^+(X;\pi)$
for all $t \geq 0$.
If $X$ is continuous, these two quantities coincide.

\begin{definition} \label{2d}
Let $E$ be a Banach space, $X \in D(\mathbb{R}_{\geq 0},E)$, and $(\pi_n)_{n \in \mathbb{N}}$
be a sequence of partitions of $\mathbb{R}_{\geq 0}$.
\begin{enumerate}
    \item The sequence $(\pi_n)$ \emph{satisfies Condition (C) for $X$}
        if it satisfies (C1)--(C3) below.
        \begin{enumerate}
            \item[(C1)] Let $t \in [0,\infty[$ and $\varepsilon > 0$. 
                Then there exists an $N \in \mathbb{N}$ such that 
                for all $n \geq N$ and for all $I \in \pi_n$,
                the set $I \cap [0,t] \cap D_{\varepsilon}(X)$
                has at most one element.
            \item[(C2)] Let $s \in D(X)$ and $t \in [s,\infty \rlbrack$. Then
                \begin{equation*}
                    \lim_{n \to \infty} \delta X_t(\pi_n(s)) 
                    =
                    \lim_{n \to \infty} \left\{ X(\overline{\pi_n}(s) \wedge t) - X(\underline{\pi_n}(s) \wedge t) \right\}
                    =
                    \Delta X_s.
                \end{equation*}
            \item[(C3)] For all $t \in \mathbb{R}_{\geq 0}$,
                \begin{equation*}
                \varlimsup_{\varepsilon \downarrow\downarrow 0} \varlimsup_{n \to \infty} O^{+}_t(X-J_{\varepsilon}(X);\pi_n) = 0.
                \end{equation*}
            \end{enumerate}
    \item The sequence \emph{$(\pi_n)$ approximates $X \colon \mathbb{R}_{\geq 0} \to E$ from the left}
        if $\lim_{n \to \infty} X(\underline{\pi_n}(t)) = X(t-)$ holds for all $t > 0$.
        Then we call $(\pi_n)$ a \emph{left approximation sequence}
        for $X$.
\end{enumerate}
\end{definition}

\begin{remark} \label{2db}
\begin{enumerate}
    \item Let $(\pi_n)$ be a sequence of partitions along which
        we consider a quadratic variation of a path.
        In this paper, we often require that $(\pi_n)$ satisfies Condition (C) defined above.
        Therefore, we can say that (C) is a condition for
        integrators of It\^{o}--F\"{o}llmer integration.
    \item In contrast to (i), Condition (ii) of Definition~\ref{2d} 
        needs to be satisfied by integrands of It\^{o}--F\"{o}llmer or Stieltjes integrals.
        This will be mainly used in Theorem~\ref{2e} and Lemma~\ref{2g}. 
\end{enumerate}
\end{remark}

Under these assumptions, we have the following $C^{1,2}$ type It{\^o}
formula for Banach space-valued paths.

\begin{theorem}[It{\^o} formula] \label{2e}
Let $E$, $E_1$, $F$, and $G$ be Banach spaces,
$B \colon E \times E \to E_1$ be a bounded bilinear map,
$(A,X) \in D(\mathbb{R}_{\geq 0}, F \times E)$,
and $(\pi_n)$ be a sequence of partitions that satisfies Condition (C) for $(A,X)$
and approximates $(A,X)$ from the left.
Suppose that $X$ has weak $B$-quadratic variation
and finite 2-variation along $(\pi_n)$
and suppose also that $A$ has finite variation.

Moreover, let $f \colon F \times E \to G$ be a function satisfying the following conditions:
\begin{enumerate}
    \item the map $F \ni a \mapsto f(a,x) \in G$ is G\^{a}teaux differentiable for all $x \in E$
        and $D_a f \colon F \times E \to \mathcal{L}_{\mathrm{c}}(F,G)$ is continuous;
    \item the map $E \ni x \mapsto f(a,x) \in G$ is twice G\^{a}teaux differentiable for all $a \in F$
        and $D_x f \colon F \times E \to \mathcal{L}_{\mathrm{c}}(E,G)$ is continuous;
    \item there is a continuous function 
        $D_B^2 f \colon F \times E \to \mathcal{L}_{\mathrm{c}}(E_1,G)$ 
        that commutes the diagram
        \begin{equation*}
            \begin{tikzcd}[column sep=42pt]
                F \times E  \arrow[r,"{D_x^2 f}"] \arrow[rd,"{D_B^2 f}"']
                & \mathcal{L}_{\mathrm{c}}(E \widehat{\otimes} E,G)   \\
                & \mathcal{L}_{\mathrm{c}}(E_1,G)  \arrow[u,"B^*"']
            \end{tikzcd}
        \end{equation*}
        where $B^* \colon \mathcal{L}_{\mathrm{c}}(E_1,G) \to \mathcal{L}_{\mathrm{c}}(E \widehat{\otimes} E,G)$
        is defined by $B^*(T)(x \otimes y) = T \circ B(x,y)$.
\end{enumerate}
Then the It\^{o}--F\"{o}llmer integral 
$\int_0^t D_x f (A_{s-},X_{s-}) \mathrm{d}X_s$ exists in the weak topology
and it satisfies
\begin{multline} \label{2f}
    f(A_t,X_t) - f(A_0,X_0)
    = 
    \int_{0}^{t} D_a f (A_{s-},X_{s-}) \mathrm{d}A^{\mathrm{c}}_s
    + \int_0^t D_x f(A_{s-},X_{s-})\mathrm{d}X_s   \\
    + \frac{1}{2} \int_{0}^{t} D_B^2 f(A_{s-},X_{s-})\mathrm{d}Q_B(X,X)^{\mathrm{c}}_s
    + \sum_{0 < s \leq t} \left\{ \Delta f(A_s,X_s) - D_x f(A_{s-},X_{s-})\Delta X_s \right\}.
\end{multline}

Moreover, if the quadratic variation $Q_B$ exits in the strong sense,
the convergence of the It\^{o}--F\"{o}llmer integral
holds in the norm topology of $G$.
\end{theorem}

Here, note that \eqref{2f} is an equation in the Banach space $G$.
The It\^{o}--F\"{o}llmer integral in Theorem~\ref{2e} is defined by 
\begin{equation*}
    \int_0^t D_x f(A_{s-},X_{s-})\,\mathrm{d}X_s
    =
    \lim_{n \to \infty} \sum_{\lrbrack r,s \rbrack \in \pi_n}
        D_x f(A_{r},X_{r})(X_{t \wedge s} - X_{t \wedge r})
\end{equation*}
with suitable topology (see Definition~\ref{5.9b}.)

As a direct consequence of Theorem~\ref{2e}, we can derive the It\^{o} formula
related to tensor quadratic variations.

\begin{corollary} \label{2fab}
Let $E$ and $F$ be Banach spaces and let $\alpha$ be a reasonable cross norm on $E \otimes E$.
Suppose that $X \in D(\mathbb{R}_{\geq 0},E)$ has strong or weak $\alpha$-tensor quadratic variation
and finite 2-variation along $(\pi_n)$
and suppose also that $A \in FV(\mathbb{R}_{\geq 0},F)$.
Moreover, let $f \colon F \times E \to G$ be a function satisfying the following conditions:
\begin{enumerate}
    \item the map $F \ni a \mapsto f(a,x) \in G$ is G\^{a}teaux differentiable for all $x \in E$
        and $D_a f \colon F \times E \to \mathcal{L}_{\mathrm{c}}(F,G)$ is continuous;
    \item the map $E \ni x \mapsto f(a,x) \in G$ is twice G\^{a}teaux differentiable for all $a \in F$
        and $D_x f \colon F \times E \to \mathcal{L}_{\mathrm{c}}(E,G)$ is continuous;
    \item the second derivative of $f$ induces a continuous function
        $D_x^2 f \colon F \times E \to \mathcal{L}_{\mathrm{c}} (E \widehat{\otimes}_{\alpha} E,G)$.
\end{enumerate}
Then $f(A,X)$ admits the following It\^o formula:
\begin{multline*}
    f(A_t,X_t) - f(A_0,X_0)
    = 
    \int_{0}^{t} D_a f (A_{s-},X_{s-}) \mathrm{d}A^{\mathrm{c}}_s
    + \int_0^t D_x f(A_{s-},X_{s-})\mathrm{d}X_s   \\
    + \frac{1}{2} \int_{0}^{t} D_x^2 f(A_{s-},X_{s-})\mathrm{d}{^\alpha [X,X]}^{\mathrm{c}}_s
    + \sum_{0 < s \leq t} \left\{ \Delta f(A_s,X_s) - D_x f(A_{s-},X_{s-})\Delta X_s \right\}.
\end{multline*}
The convergence of the It\^{o}--F\"{o}llmer integral holds in the norm or weak topology, respectively.
\end{corollary}

\begin{remark} \label{2fac}
\begin{enumerate}
    \item Let $f \colon F \times E \to G$ be a $C^{1,2}$ function in the sense of Fr\'{e}chet differentiation,
        i.e.
        \begin{enumerate}
            \item the function $a \mapsto f(a,x)$ is Fr\'{e}chet differentiable for all $x \in E$ 
                and $D_a f \colon F \times E \to \mathcal{L}(F,G)$ is continuous;
            \item the function $x \mapsto f(a,x)$ is twice Fr\'{e}chet differentiable for all $a \in F$,
                and $Df \colon F \times E \to \mathcal{L}(E,G)$ and 
                $D^2 f \colon F \times E \to \mathcal{L}(E \widehat{\otimes} E,G)$ are continuous.
        \end{enumerate}
        Then $f$ satisfies the assumptions in Corollary~\ref{2fab}
        for the projective tensor norm $\gamma$.
    \item Define a class $FC^2_{\mathrm{b}}(E)$ to be the set of all real functions of the form
        $G = g \circ P$ with $P \colon E \to V$ being the projection to a finite-dimensional subspace 
        and $g \in C^2_{\mathrm{b}}(V)$. Let $\mathcal{C}$ be the completion of 
        $FC^2_{\mathrm{b}}$ with respect to the norm
            \begin{equation*}
                \lVert G \rVert
                \coloneqq
                \sup_{x \in E} \lvert G(x) \rvert 
                + \sup_{x \in E} \lVert DG(x) \rVert_{E^*} 
                + \sup_{x \in E} \lVert D^2 G(x) \rVert_{(E \widehat{\otimes}_{\alpha} E)^*}.
            \end{equation*}
        Then every element of $\mathcal{C}$ satisfies the assumptions in Corollary~\ref{2fab}.
\end{enumerate}
\end{remark}

\begin{remark} \label{2fb}
A possible direction of development is a functional extension of
Theorem~\ref{2e} and Corollary~\ref{2fab},
that is, infinite-dimensional functional It\^{o} calculus.
For such an extension, we should introduce Dupire's derivative of
a functional of infinite dimensional paths.
This seems an important future work.
See, as referred to in Section 1,
Dupire~\cite{Dupire_2009}, Cont and Fourni{\'e}~\cite{Cont_Fournie_2010,
Cont_Fournie_2010b,
Cont_Fournie_2013}, and 
Ananova and Cont~\cite{Ananova_Cont_2017}, for 
functional It\^{o} calculus in finite-dimensional state space.
\end{remark}

The following lemma is essentially used to prove Theorem~\ref{2e}.

\begin{lemma} \label{2g}
Let $B \colon E \times E \to E_1$ be a bounded bilinear map between Banach spaces
and $(\pi_n)$ be a sequence of partitions of $\mathbb{R}_{\geq 0}$.
Suppose that a path $X \in D(\mathbb{R}_{\geq 0},E)$ has
weak $B$-quadratic variation and finite $2$-variation along $(\pi_n)$.
Moreover, assume that $(\pi_n)$ 
satisfies Condition (C) for $X$ and approximates 
a $\xi \in D(\mathbb{R}_{\geq 0},\mathcal{L}_{\mathrm{c}}(E_1,G))$
from the left.
Then for all $t \in [0,\infty[$
\begin{equation} \label{2h}
    \sum_{I \in \pi_n} 
        \xi B(\delta X_t,\delta X_t)(I)
    \xrightarrow[n \to \infty]{}
    \int_{\lrbrack 0,t \rbrack} \xi_{u-} \mathrm{d}Q_B(X,X)_u
\end{equation}
holds in the weak topology of $G$.

If, in addition, $Q_B(X,X)$ is the strong $B$-quadratic variation, 
then \eqref{2h} holds in the norm topology.
\end{lemma}

\section{Auxiliary results regarding sequences of partitions}

In this section, we investigate conditions on a sequence of partitions 
along which we deal with quadratic variations and It\^{o}--F\"{o}llmer integrals.
Recall that basic notions were defined in Definition~\ref{2d}.

\begin{definition} \label{4b}
Let $E$ be a Banach space, $X \in D(\mathbb{R}_{\geq 0},E)$, and $(\pi_n)_{n \in \mathbb{N}}$
be a sequence of partitions of $\mathbb{R}_{\geq 0}$.
\begin{enumerate}
\item The sequence $(\pi_n)$
    \emph{controls the oscillation of $X$}
    if $\lim_{n \to \infty}O^-_t(X;\pi_n) = 0$ holds for all $t$.
\item The sequence $(\pi_n)$ \emph{exhausts the jumps of $X$} if
    $D(X) \subset \bigcup_{n \in \mathbb{N}} \bigcap_{k \geq n} \pi^{\mathrm{p}}_k$.
\end{enumerate}
\end{definition}

\begin{example} \label{4c}
\begin{enumerate}
    \item Let $r$ be an irrational number and $X = 1_{\lbrack r,\infty \rlbrack}$.
        For each $n \in \mathbb{N}$, we set $\pi_n = \{ \lrbrack k2^{-n},(k+1)2^{-n} \rbrack ; k \in \mathbb{N} \}$.
        Then the sequence $(\pi_n)$ satisfies $\lvert \pi_n \rvert \to 0$ as $n \to \infty$.
        This sequence, however, does not control the oscillation of $X$.
    \item Let $X = 1_{\lbrack 1,\infty \rlbrack}$
        and $\pi_n = \{ \lrbrack k,k+1 \rbrack; k \in \mathbb{N} \}$.
        Though the sequence $(\pi_n)$ controls the oscillation of $X$,
        it does not satisfy $\lvert \pi_n \rvert \to 0$.
    \item Let $X = 1_{\lbrack 1/2, \infty \rlbrack}$
        and $\pi_n = \{ \lrbrack k,k+1 \rbrack; k \in \mathbb{N} \}$.
        The sequence $(\pi_n)$ neither controls the oscillation of $X$
        nor satisfies $\lvert \pi_n \rvert \to 0$.
        However, it satisfies Condition (C) for $X$ in the sense of Definition~\ref{2d}.
\end{enumerate}
\end{example}

Condition~(i) of Definition~\ref{4b} is characterized as follows.

\begin{lemma} \label{4e}
Let $X$ be a c\`{a}dl\`{a}g path in $E$ and $(\pi_n)$ be a sequence of partitions of $\mathbb{R}_{\geq 0}$.
Then the following conditions are equivalent.
\begin{enumerate}
    \item The sequence $(\pi_n)$ controls the oscillation of $X$.
    \item The sequence $(\pi_n)$ satisfies the following conditions:
        \begin{enumerate}
            \item the sequence $(\pi_n)$ exhausts the jumps of $X$;
            \item if $X$ is not constant on $\mathopen{]}s,t\mathclose{[} \subset \mathbb{R}_{\geq 0}$, 
                then $\lrbrack s,t \rlbrack$ contains at least one element of
                $\pi_n^{\mathrm{p}}$ for sufficiently large $n$.
        \end{enumerate}
\end{enumerate}
\end{lemma}

Lemma~\ref{4e} is a generalization of
Lemma~1 of Vovk~\cite[272]{Vovk_2015}.
Note that the condition
`\emph{$]s,t[$ contains at least one point of $\pi_n^{\mathrm{p}}$}'
is equivalent to
`\emph{there is no $I \in \pi_n$ including $\lrbrack s,t \rlbrack$}.'

\begin{proof}
\emph{Step~1.1: (i) $\Longrightarrow$ (ii)-(a).}
Let $s \in D(X)$ and set $\varepsilon = \lVert \Delta X(s) \rVert_{E}$.
Moreover, fix $T > s$ arbitrarily.
Then, by assumption, we can choose an $N \in \mathbb{N}$ such that
$O_T^-(X,\pi_n) < \varepsilon/2$ holds for all $n \geq N$.
We will show that $s = \overline{\pi_n}(s)$ holds for all $n \geq N$.
Take an $s'$ from the interval $\mathopen{]}\underline{\pi_n}(s),s \mathclose{[}$
such that
\begin{equation*}
    \lVert X_{s'} - X_{s-} \rVert_{E} \leq O^{-}_T(X,\pi_n) < \frac{\varepsilon}{2}.
\end{equation*}
Then, we have
\begin{equation*}
    \lVert X_s - X_{s'} \rVert_{E} 
    \geq 
    \lVert \Delta X_s \rVert_{E} - \lVert X_{s-} - X_{s'} \rVert_{E} 
    > \varepsilon - \frac{\varepsilon}{2} 
    = 
    \frac{\varepsilon}{2}.
\end{equation*}
This shows that $s \notin \lrbrack \underline{\pi_n}(s),\overline{\pi_n}(s) \rlbrack$,
and therefore $s = \overline{\pi_n}(s)$.
Hence $(\pi_n)$ exhausts the jumps of $X$.

\emph{Step~1.2: (i) $\Longrightarrow$ (ii)-(b).}
Assume that $\varepsilon \coloneqq \omega(X;\mathopen{]} s,t \mathclose{[}) > 0$.
Choose an $N \in \mathbb{N}$ that satisfies 
$O^{-}_t(X,\pi_n) < \varepsilon$ for all $n \geq N$.
For an arbitrarily fixed $n \geq N$, choose a unique $i$ satisfying $s \in [t^n_i,t^n_{i+1}[$.
It remains to show that $t^n_{i+1} \in \mathopen{]}s,t \mathclose{[}$.
Seeking a contradiction, suppose $t^n_{i+1} \geq t$.
Then $\lbrack s,t \rlbrack \subset \lbrack t^n_i,t^n_{i+1} \rlbrack$ and hence
\begin{equation*}
    O^{-}_t(X,\pi_n) 
    \geq
    \sup_{u,v \in [t^n_i,t^n_{i+1}\mathclose{[} \cap [0,t]} \lVert X_u - X_v \rVert
    \geq
    \sup_{u,v \in [s,t[} \lVert X_u - X_v \rVert
    =
    \varepsilon.
\end{equation*}
This contradicts the condition that $O^{-}_t(X,\pi_n) < \varepsilon$.

\emph{Step~2: (ii) $\Longrightarrow$ (i).}
Suppose that $(\pi_n)$ satisfies the conditions (ii)-(a) and (b).

Fix an $\varepsilon > 0$ and a $t > 0$ arbitrarily.
Because $X$ is c\`{a}dl\`{a}g, we can take a sequence $0 = s_0 < s_1 < \dots < s_N = t$
such that $\omega(X; \mathopen{]} s_i,s_{i+1}\mathclose{[}) < \varepsilon/2$
for all $i$ (see Lemma~\ref{3b}).
By assumption, we can choose an $N \in \mathbb{N}$ satisfying
the following conditions:
\begin{enumerate}
    \item[1.] If $n \geq N$, there are no $I \in \pi_n$ and $i \in \{ 0,\dots, N \}$
        satisfying $\lrbrack s_i,s_{i+1} \rlbrack \subset I$ and
        $\omega(X;\lrbrack s_i,s_{i+1} \rlbrack) > 0$.
    \item[2.] $
        \{ s_0,\dots, s_N \} \cap D(X) 
        \subset \bigcap_{n \geq N} \pi_n^{\mathrm{p}}
        $.
\end{enumerate}

Let $n \geq N$ and $\lrbrack u,v \rbrack \in \pi_n$.
First, assume that $\omega(X; \lrbrack u,v \rlbrack) > 0$.
By Condition~1, we see that there are only two cases
for the relationship between $\lrbrack u,v \rbrack$ and $(s_i)_{0 \leq i \leq N}$ as follows.
\begin{enumerate}
    \item[A.] There is a unique $i$ such that $\lrbrack u,v \rbrack \subset \lrbrack s_i,s_{i+1}\rlbrack$.
    \item[B.] There is a unique $i$ such that $s_i \in \lrbrack u,v \rlbrack$.
\end{enumerate}
In Case~A, the oscillation of $X$ on $\lrbrack u,v \rlbrack$ is estimated as
\begin{equation*}
    \omega(X;\lrbrack u,v \rlbrack \cap [0,t]) \leq \omega(X; \lrbrack s_i,s_{i+1}\rlbrack) < \frac{\varepsilon}{2}.
\end{equation*}
On the other hand, in Case~B, 
$X$ is continuous at $s_i \in \lrbrack u,v \rlbrack$ by Condition~2. Therefore,
\begin{align*}
    \omega(X;\lrbrack u,v \rlbrack \cap [0,t])
    \leq
    \omega(X; \lrbrack s_{i-1},s_{i}\rlbrack) + \omega(X; \lrbrack s_i,s_{i+1}\rlbrack) 
    <
    \varepsilon.
\end{align*}
If $\omega(X; \lrbrack u,v \rlbrack) = 0$, we clearly have the same estimate.
By the discussion above,
we find that $\omega(X;\lrbrack u,v \rlbrack \cap [0,t]) < \varepsilon$ holds for all
$\lrbrack u,v \rbrack \in \pi_n$, and consequently
\begin{equation*}
    O^{-}_t(X;\pi_n) = \sup_{\lrbrack r,s \rbrack \in \pi_n} \omega(X;\lrbrack r,s \rlbrack \cap [0,t]) \leq \varepsilon
\end{equation*}
for every $n \geq N$. This completes the proof.
\end{proof}

We next consider the condition that $(\pi_n)$ approximates a path from the left.
Given a partition $\pi$, define the left discretization of
a path $\xi \colon \mathbb{R}_{\geq 0} \to E$ along $\pi$ by
\begin{equation*}
    {^\pi\xi} = \sum_{\mathopen{]}r,s] \in \pi} \xi(r) 1_{]r,s]}.
\end{equation*}
If $\xi$ is c\`{a}dl\`{a}g, then the sequence $(\pi_n)$ approximates $\xi$ from the left
in the sense of Definition~\ref{2d}
if and only if $({^{\pi_n}\xi})_{n \in \mathbb{N}}$ converges
to $\xi_{-}$ pointwise.
Consider the following example.

\begin{example} \label{4d}
Let $X = 1_{\lbrack 1/2, \infty \rlbrack}$
and $\pi_n = \{ \lrbrack k,k+1 \rbrack; k \in \mathbb{N} \}$.
As we saw in Example~\ref{4c},
the sequence $(\pi_n)$ satisfies Condition (C) for $X$.
For each $n \in \mathbb{N}$,
the left discretization of $X$ is given by 
${^{\pi_n} X} = 1_{\lrbrack 1,\infty \rlbrack}$.
The sequence $({^{\pi_n} X})$ does not converge to $X$ pointwise,
and hence $(\pi_n)$ does not approximate $X$ from the left.
\end{example}

As we mentioned in Section~1, 
two types of assumptions about a sequence of partitions are 
frequently used in the context of the It\^{o}--F\"{o}llmer calculus.
One is that $\lvert \pi_n \rvert \to 0$ and the other is that $(\pi_n)$
controls the oscillation of $X$.
In the next proposition, we show that both conditions 
imply that $(\pi_n)$ satisfies Condition (C) for $X$.

\begin{proposition} \label{4f}
Let $(\pi_n)$ be a sequence of partitions of $\mathbb{R}_{\geq 0}$ and let $E$ be a Banach space.
\begin{enumerate}
    \item Suppose that $(\pi_n)$ satisfies $\lvert \pi_n \rvert \to 0$. 
        Then it satisfies Condition (C) for every c\`{a}dl\`{a}g path in $E$.
        Moreover, it approximates every c\`{a}dl\`{a}g path in $E$ from the left.
    \item Suppose that $(\pi_n)$ controls the oscillation of $X \in D(\mathbb{R}_{\geq 0},E)$. 
        Then it satisfies Condition (C) for $X$ and approximates $X$ from the left.
\end{enumerate}
\end{proposition}

\begin{proof}
(i)
If $\lvert \pi_n \rvert \to 0$, then $\overline{\pi_n}(t) \to t$
and $\underline{\pi_n}(t) \to t$ hold for every $t \geq 0$.
This directly implies that $(\pi_n)$ approximates $X$ from the left.
Moreover, we have $\delta X_u(\pi_n(t)) \to X_{t}$ for every $t,u > 0$
with $t \leq u$. Hence, $(\pi_n)$ satisfies Condition~(C2).
Condition~(C3) follows from (iii) of Lemma~\ref{3b}.
Condition~(C1) remains to be shown.
Given an $\varepsilon > 0$, define
\begin{equation*}
    r \coloneqq \inf \{ \lvert u-v \rvert \mid u,v \in D_{\varepsilon}(X) \cap [0,t], u \neq v \} > 0.
\end{equation*}
If $D_{\varepsilon}(X) \cap [0,t]$ has only one element,
there is nothing to do.
Otherwise, $r$ is not zero, 
because $D_{\varepsilon}(X) \cap [0,t]$ has at most finitely many elements
(see (i) of Lemma~\ref{3b}).
Now we take an $N$ satisfying $\lvert \pi_n \rvert < r$ for all $n \geq N$.
Then for each $n \geq N$ and $\lrbrack u,v \rbrack \in \pi_n$, 
the set $\lrbrack u,v \rbrack \cap [0,t]$ contains at most one element of $D_{\varepsilon}$.

(ii)
Assume that $(\pi_n)$ controls the oscillation of $X$.
Then we see that $(\pi_n)$ approximates $X$ from the left
by the following estimate.
\begin{equation*}
    \lVert X_{\underline{\pi_n}(t)} - X_{t-} \rVert_{E}
    \leq
    \omega(X;\lbrack \underline{\pi_n}(t), \overline{\pi_n}(t)\rlbrack \cap [0,t] )
    \leq
    O^{-}_t(X,\pi_n).
\end{equation*}

Now, let us show that $(\pi_n)$ satisfies (C) for $X$.
To obtain (C2), take a $t \in D(X)$.
Because $(\pi_n)$ exhausts the jumps of $X$ (Lemma~\ref{4e}), 
we have $\overline{\pi_n}(t) = t$ for sufficiently large $n$.
This combined with the fact that $(\pi_n)$ is a left-approximation sequence
implies (C2).
Next, consider (C1). 
Let $\varepsilon > 0$ and fix an $N_{\varepsilon} \in \mathbb{N}$
so that $O^{-}_{t}(X,\pi_n) < \varepsilon$ holds for any $n \geq N_{\varepsilon}$.
Then, for every $n \geq N_{\varepsilon}$
and $\lrbrack r,s] \in \pi_n$, the interval $\lrbrack r,s \rlbrack \cap [0,t]$ does not 
contain any jump of $X$ that is greater than $\varepsilon$.
Therefore, $\lrbrack r,s \rbrack$ possesses at most one element of $D_{\varepsilon}(X)$.
This means that $(\pi_n)$ satisfies (C1).
All that is left is to check Condition~(C3).
Choose an $M_{\varepsilon}$
that satisfies $O^{-}_t(X;\pi_n) < \varepsilon/2$ for all $I \in \pi_n$ and $n \geq M$.
As we just have shown, $J_{\varepsilon/2}(X)$ is zero on the interior of each 
$I \in \pi_n$ and $n \geq M_{\varepsilon}$. Hence, 
\begin{align*}
    \omega(X-J_{\varepsilon/2}(X);\lrbrack r,s \rbrack \cap [0,t])
& \leq
    \omega(X-J_{\varepsilon/2}(X);\lrbrack r,s \rlbrack \cap [0,t]) + \lVert \Delta (X-J_{\varepsilon/2}(X))_s \rVert_{E}  \\
& \leq 
    O^{-}_t(X;\pi_n)
    + \lVert \Delta (X-J_{\varepsilon/2}(X))_s \rVert_{E}
<
    \varepsilon
\end{align*}
holds for all $\lrbrack r,s \rbrack \in \pi_n$ and $n \geq M_{\varepsilon}$,
which implies (C3).
\end{proof}

In the last part of this section,
we give an additional lemma about a sequence of partitions.

\begin{lemma} \label{4g}
\begin{enumerate}
    \item Let $X$ be a c\`{a}dl\`{a}g path in a Banach space $E$.
        If $(\pi_n)$ approximates $\xi$ from the left,
        then $(\pi_n)$ also approximates $f \circ \xi$ from the left for every
        continuous function $f \colon E \to E'$ to an arbitrary Banach space.
    \item Let $X$ and $Y$ be c\`{a}dl\`{a}g paths in Banach spaces $E$ and $F$,
        respectively. If $(\pi_n)$ satisfies Condition (C) for the path $(X,Y)$ in $E \times F$,
        then $(\pi_n)$ satisfies (C) for each of $X$ and $Y$.
        Here, we regard $E \times F$ as a Banach space endowed with the direct sum norm
        $\lVert \phantom{x} \rVert_E + \lVert \phantom{x} \rVert_F$.
\end{enumerate}
\end{lemma}

\begin{proof}
(i) immediately follows from the continuity of $f$.

To show (ii), suppose that $(\pi_n)$ satisfies (C) for $(X,Y)$.
It suffices to show that $(\pi_n)$ satisfies (C) for $X$.
First, fix $t \in \mathbb{R}_{\geq 0}$ and $\varepsilon > 0$ arbitrarily,
and then choose an $N$ so that $D_\varepsilon(X,Y) \cap I \cap [0,t]$
has at most one element for all $n \geq N$ and $I \in \pi_n$.
The inclusion
$D_\varepsilon(X) \cap I \cap [0,t] \subset D_\varepsilon(X,Y) \cap I \cap [0,t]$
implies that the cardinality of $I \cap [0,t] \cap D_\varepsilon(X)$ is 
no greater than 1.
Condition~(C2) follows directly from the definition of product topology.
Condition~(C3) remains to be shown.
For an arbitrary $\delta > 0$, 
choose an $\varepsilon_0 > 0$ satisfying
\begin{equation*}
    \sup_{\varepsilon \leq \varepsilon_0} \varlimsup_{n \to \infty} O^+_t((X,Y)-J_\varepsilon(X,Y);\pi_n) < \frac{\delta}{2}.
\end{equation*}
Set $\varepsilon_1 = \varepsilon_0 \wedge (\delta/2)$.
Given $\varepsilon \leq \varepsilon_1$,
we can take $M_{\varepsilon} > 0$ such that
\begin{enumerate}
    \item[(a)] $I \cap [0,t] \cap D_{\varepsilon}(X,Y)$ has at most one element
        for all $I \in \pi_n$ and $n \geq M_{\varepsilon}$.
    \item[(b)] $\sup_{n \geq M_{\varepsilon}} O^+_t((X,Y)-J_{\varepsilon}(X,Y);\pi_n) < \delta/2$.
\end{enumerate}
If $n \geq M_{\varepsilon}$ and $I \in \pi_n$, then for any $u,v \in I$, we have
\begin{align*}
    & 
    \lVert (X-J_\varepsilon(X))_u - (X-J_\varepsilon(X))_v \rVert_{E}  \\
    & \qquad \leq 
    \lVert ( X-J(D_\varepsilon(X,Y),X) )_u - ( X-J (D_\varepsilon(X,Y),X) )_v \rVert_{E}  \\
    & \qquad\quad 
    + \lVert ( J(D_\varepsilon(X,Y),X) - J(D_\varepsilon(X),X) )_u - ( J(D_\varepsilon(X,Y),X) - J (D_\varepsilon(X),X) )_v \rVert_{E}  \\
    & \qquad \leq
    \sup_{n \geq M_{\varepsilon}} O^+_t((X,Y)-J_{\varepsilon}(X,Y);\pi_n) + \varepsilon   \\
    & \qquad \leq
    \frac{\delta}{2} + \frac{\delta}{2} = \delta.
\end{align*}
Here, note that the second inequality holds by Condition~(a) above.
Thus, we get
\begin{equation*}
    \varlimsup_{n \to \infty} O^+_t(X-J_{\varepsilon}(X);\pi_n)
    \leq
    \sup_{n \geq M_{\varepsilon}} O^+_t(X-J_{\varepsilon}(X);\pi_n)
    \leq
    \delta.
\end{equation*}
for arbitrary $\varepsilon \leq \varepsilon_1$.
This implies (C3) for $X$.
\end{proof}

\section{Properties of quadratic variations}

This section is devoted to studying some basic properties
of quadratic variations introduced in Section 3.
Throughout this section, suppose that
we are given a sequence $\Pi = (\pi_n)_{n \in \mathbb{N}}$
of partitions of $\mathbb{R}_{\geq 0}$.

First, we give some examples of quadratic variations.

\begin{example} \label{5b}
Let $A$ be a path of finite variation in a Banach space $E$.
If $(\pi_n)$ satisfies $\lvert \pi_n \rvert \to 0$,
then $A$ has projective tensor and scalar quadratic variations given by 
\begin{equation*}
    {^\gamma [A,A]_t} = \sum_{0 < s \leq t} (\Delta A_s)^{\otimes 2}, \qquad 
    Q(A)_t = \sum_{0 < s \leq t} \lVert \Delta A_s \rVert^2.
\end{equation*}
This result will be proved later in this section.
\end{example}

\begin{example} \label{5c}
Let $x \colon [0,\infty \rlbrack \to \mathbb{R}$
be a c\`{a}dl\`{a}g path and $(\pi_n)$ a sequence of partitions 
such that $\lvert \pi_n \rvert \to 0$.
Given arbitrary $C^1$-function $f \colon \mathbb{R} \to E$
into a Banach space $E$, let us set $X(t) = f(x(t))$.
If $x$ has quadratic variation along $(\pi_n)$,
then $X$ has projective tensor quadratic variation given by
\begin{equation*}
    {^\gamma [X,X]}_t = \int_0^t Df(x_{s-})^{\otimes 2} \mathrm{d}[x,x]^{\mathrm{c}}_s
    + \sum_{0 < s \leq t} \Delta f(x_s)^{\otimes 2}.
\end{equation*}
This is one of the main results of the second article 
in this series~\cite{Hirai_2022b}.

For the construction of 
a real continuous path of nontrivial quadratic variation,
refer to Schied~\cite{Schied_2016}, 
Mishura and Schied~\cite{Mishura_Schied_2016},
and Cont and Das~\cite{Cont_Das_2022}.
\end{example}

The next examples are from the theory of stochastic processes.

\begin{example} \label{5d}
\begin{enumerate}
    \item Let $(\Omega,\mathcal{F},(\mathcal{F}_t)_{t \geq 0},P)$ be a filtered probability space satisfying the usual conditions.
        Consider a semimartingale $X = (X_t)_{t \geq 0}$ 
        in a separable Hilbert space $H$.
        Moreover, let $\pi = (\tau^n_k)_{n \in \mathbb{N}}$ be an 
        increasing sequence of bounded stopping times such that 
        $\tau^n_k \to \infty$ as $k \to \infty$
        and $\sup_{k} (\tau^{n+1}_k - \tau^n_k) \to 0$ as $n \to \infty$
        almost surely.
        Then the process $[X,X]^{\pi}$ converges to 
        the quadratic variation process $[X,X]$ uniformly on compacts in probability (ucp).
        By passing to an appropriate subsequence, 
        we see that almost all paths have quadratic variation 
        along the subsequence.
        See Gravereaux and Pellaumail~\cite{Gravereaux_Pellaumail_1974} or
        Metivier and Pellaumail~\cite{Metivier_Pellaumail_1980b} for details.
    \item In addition to the assumptions of (i), 
        let $f \colon H \to E$ be a $C^1$ function into a Banach space $E$.
        Then, along a suitable subsequence of $(\pi_n)$,
        almost all paths of $f(X)$ have quadratic variation.
        If, moreover, $f$ is of $C^2$ class, the It\^{o}--F\"{o}llmer integral 
        $\int_0^{\cdot} Df(X_{s-}) \mathrm{d}X_s$ exists and its paths have 
        quadratic variation along the same subsequence.
        Now, since the Banach space $E$ is chosen arbitrarily,
        it may fail to satisfy some useful properties
        required by the martingale theory,
        e.g. UMD property or martingale type 2 property.
        The path $f(X)$, however, behaves well enough 
        from the viewpoint of the It\^{o}--F\"{o}llmer calculus. 
\end{enumerate}
\end{example}

Now we consider the transpose of quadratic covariation.
Let $B \in \mathcal{L}^{(2)}(E,F;G)$ and define the transpose $^tB \colon F \times E \to G$
of $B$ by $^tB(y,x) = B(x,y)$.
Then $(X,Y) \in D(\mathbb{R}_{\geq 0},E \times F)$
has strong/weak $B$-quadratic covariation if and only if $(Y,X)$ does
with respect to the transpose $^tB$.

Recall that a $d$-dimensional c\`{a}dl\`{a}g path 
$X = (X_1,\dots, X_d)$ has tensor quadratic variation 
if and only if $X_i$ and $X_j$ have quadratic covariation 
for each $i$ and $j$.
This characterization is generalized to bilinear quadratic covariations in Banach spaces.

\begin{proposition} \label{5.4c}
Let $E_i$, $F_j$, and $G_{ij}$ be Banach spaces
and let $B_{ij} \colon E_i \times F_j \to G_{ij}$ be a bounded bilinear map
for $i,j \in \{ 1,2\}$.
Define a continuous bilinear map $\mathbf{B} \colon (E_1 \times E_2) 
\times (F_1 \times F_2) \to \prod_{i,j} G_{ij}$ by
\begin{equation*}
    \mathbf{B}((x_1,x_2),(y_1,y_2)) = (B_{ij}(x_i,y_i))_{i,j \in \{1,2 \}}.
\end{equation*}
If $\mathbf{X} = (X_1,X_2) \in D(\mathbb{R}_{\geq 0}, E_1 \times E_2)$
and $\mathbf{Y} = (Y_1,Y_2) \in D(\mathbb{R}_{\geq 0}, F_1 \times F_2)$,
then $(\mathbf{X},\mathbf{Y})$ has strong or weak 
$\mathbf{B}$-quadratic covariation
if and only if $(X_i,Y_j)$ has respectively
strong or weak $B_{ij}$-quadratic covariation for all $i,j \in \{1,2\}$.
In this case, these quadratic covariations satisfy
\begin{equation} \label{5.4cb}
    Q_{\mathbf{B}}((X_1,X_2),(Y_1,Y_2)) = (Q_{B_{ij}}(X_i,X_j))_{i,j \in \{1, 2 \}}.
\end{equation}
\end{proposition}

Using the matrix notation, we can also express equation~\eqref{5.4cb} as
\begin{equation*}
    Q_{\mathbf{B}}((X_1,X_2),(Y_1,Y_2)) =
    \begin{pmatrix}
        Q_{B_{11}}(X_1,Y_1) & Q_{B_{12}}(X_1,Y_2)  \\
        Q_{B_{21}}(X_2,Y_1) & Q_{B_{22}}(X_2,Y_2)
    \end{pmatrix}.
\end{equation*}

\begin{proof}
By the definition of $\mathbf{B}$, we can easily check that
\begin{gather*}
    % 1 
    \sum_{I \in \pi_n} 
        \mathbf{B} \Bigl( \bigl( (\delta X_1)_t, (\delta X_2)_t \bigr), \bigl( (\delta Y_1)_t, (\delta Y_2)_t \bigr) \Bigr)(I)
    = 
    \left( \sum_{I \in \pi_n} 
        B_{ij} \bigl( (\delta X_i)_t,(\delta Y_i)_t \bigr)(I) \right)_{i,j \in \{1,2 \}},  \\
    % 2
    \mathbf{B} \left( \Delta (X_1,X_2)_t,\Delta (Y_1,Y_2)_t \right)
    = 
    \left( B_{ij}( \Delta X_i(t),\Delta Y_j(t) ) \right)_{i,j \in \{1,2\}},
\end{gather*}
hold for all $t \geq 0$.
These immediately prove the assertion.
\end{proof}

Applying Proposition~\ref{5.4c} to
the canonical bilinear map $\otimes \colon E_i \times E_j \to E_i \widehat{\otimes}_{\alpha} E_j$,
we obtain the following corollary.

\begin{corollary} \label{5.4d}
Let $(X_1,X_2) \in D(\mathbb{R}_{\geq 0},E_1 \times E_2)$
and let $\alpha$ be a uniform crossnorm.
Then $(X_1,X_2)$ has strong or weak $\alpha$-tensor quadratic variation if and only if
$(X_i,X_j)$ has strong or weak $\alpha$-tensor quadratic covariation, respectively, 
for every $i,j \in \{1,2 \}$.
\end{corollary}

Like the quadratic covariation $[X,Y]$ of scalar paths,
quadratic covariation $Q_B$ is bilinear in an appropriate sense.

\begin{proposition} \label{5.4e}
Let $X_1,X_2 \in D(\mathbb{R}_{\geq 0},E)$ and
$Y_1,Y_2 \in D(\mathbb{R}_{\geq 0},F)$.
Suppose that $(X_i,Y_j)$ has strong or weak quadratic covariation
with respect to $B \in \mathcal{L}^{(2)} (E,F;G)$ for each $i,j \in \{ 1,2 \}$.
Then, $(X_1 + X_2,Y_1 + Y_2)$ has respectively strong or weak
$B$-quadratic covariation given by
\begin{equation*}
    Q_B(X_1+X_2,Y_1+Y_2) = \sum_{ i,j \in \{1,2\} } Q_B(X_i,Y_j).
\end{equation*}
\end{proposition}

\begin{proof}
By the bilinearity of $B$, we see that
\begin{equation*}
    \sum_{I \in \pi_n} B(\delta (X_1+X_2)_t,\delta (Y_1+Y_2)_t)(I)
    =
    \sum_{i,j \in \{1,2 \}} \sum_{I \in \pi_n} B(\delta (X_i)_t,\delta (Y_j)_t)(I)
\end{equation*}
for every $t \geq 0$. Therefore, by assumption, the left-hand side converges to $\sum_{ij} Q_B(X_i,Y_j)$
in the corresponding topology. Again by bilinearity,
\begin{equation*}
        \Delta \left( \sum_{i,j \in \{1,2\}} Q_B(X_i,Y_j) \right)(t)
    = 
        \sum_{i,j \in \{1,2\}} B(\Delta (X_i)_t, \Delta (Y_j)_t)  
    = 
        B\left( \Delta (X_1 + X_2)_t, \Delta (Y_1 + Y_2)_t \right).
\end{equation*}
Hence, $\sum_{ij} Q_B(X_i,Y_j)$ is the $B$-quadratic covariation of $X_1 + X_2$ and $Y_1 + Y_2$.
\end{proof}

\begin{corollary} \label{5.4eb}
Let $E$ and $F$ be Banach spaces and $\alpha$ be a reasonable crossnorm on $E \otimes F$.
Suppose that $(X_i,Y_j) \in D(\mathbb{R}_{\geq 0},E \times F)$ has 
strong or weak $\alpha$-tensor quadratic covariation for every $i,j \in \{ 1,2 \}$.
Then, $(X_1 + X_2,Y_1 + Y_2)$ has
strong or weak $\alpha$-tensor quadratic covariation, respectively, and it satisfies
\begin{equation*}
    {^\alpha [X_1+X_2,Y_1+Y_2]} 
    = 
    {^\alpha [X_1,Y_1]} + {^\alpha [X_1,Y_2]} + {^\alpha [X_2,Y_1]} + {^\alpha [X_2,Y_2]}.
\end{equation*}
\end{corollary}

Next, we investigate the quadratic variation of a path of finite variation.
For convenience, we introduce the following notation.
Let $D \subset \mathbb{R}_{\geq 0}$
and define functions
$e^1_D$ and $e^2_D$ from $\mathcal{P}(\mathbb{R}_{\geq 0})$ to $\{ 0,1 \}$ by
\begin{equation*}
    e^1_D(A) =
    \begin{cases}
        1 & \text{if} \quad A \cap D \neq \emptyset  \\
        0 & \text{if} \quad A \cap D = \emptyset
    \end{cases}, \qquad 
    e^2_D = 1 - e^1_D.
\end{equation*}
The symbol $\mathcal{P}(\mathbb{R}_{\geq 0})$ above denotes the 
power set of $\mathbb{R}_{\geq 0}$.

\begin{proposition} \label{5.6b}
Let $E$, $F$, and $G$ be Banach spaces and let $B \in \mathcal{L}^{(2)}(E,F;G)$.
Assume that $A \in FV(\mathbb{R}_{\geq 0},E)$,
$X \in D(\mathbb{R}_{\geq 0},F)$,
and $(\pi_n)$ satisfies Condition (C) for $(A,X)$.
Then $(A,X)$ has the strong $B$-quadratic variation given by
\begin{equation*}
    Q_B(A,X)_t = \sum_{0 < s \leq t} B(\Delta A_s,\Delta X_s).
\end{equation*}
\end{proposition}

\begin{proof}
Fix $t \in \mathbb{R}_{\geq 0}$ and 
take an arbitrary $\varepsilon > 0$.
For convenience, set $D = D(A,X)$, $D_{\varepsilon} = D_{\varepsilon}(A,X)$,
and $D^{\varepsilon} = D^{\varepsilon}(A,X)$.
Then,
\begin{align} \label{5.6c}
    & 
        \left\lVert 
            \sum_{I \in \pi_n} B(\delta A_t,\delta X_t)(I) 
            - \sum_{0 < u \leq t} B( \Delta A_u,\Delta X_u )
        \right\rVert    \notag \\
    & \qquad \leq 
        \left\lVert 
            \sum_{I \in \pi_n}
                B(\delta A_t,\delta X_t)(I) e^1_{D_{\varepsilon}}(I) 
                - \sum_{0 < u \leq t} B\left( \Delta A_u, \Delta X_u \right)
        \right\rVert
        + 
        \left\lVert 
            \sum_{I \in \pi_n} B(\delta A_t,\delta X_t)(I)e^2_{D_{ \varepsilon}}(I)
        \right\rVert
\end{align}
for any $t \in \mathbb{R}_{\geq 0}$.
We will observe the behaviour of each term on the right-hand side of \eqref{5.6c}.

Because $(\pi_n)$ satisfies Condition (C) for $X$,
there exists an $N_1$ such that $D_{\varepsilon} \cap [0,t] \cap I$
contains at most one point for all $n \geq N_1$ and $I \in \pi_n$.
If $n \geq N_1$, we have
\begin{equation*}
    \sum_{I \in \pi_n} B(\delta A_t,\delta X_t)(I) e^1_{D_{\varepsilon}}(I)
    = 
    \sum_{u \in D_{\varepsilon}} B\left( \delta A_t, \delta X_t \right)(\pi_n(u)).
\end{equation*}
Therefore, by Condition~(C2),
\begin{equation*}
    \lim_{n \to \infty}
        \sum_{I \in \pi_n} B(\delta A_t,\delta X_t)(I) e^1_{D_{\varepsilon}}(I)
    =
    \sum_{u \in D_{\varepsilon} \cap [0,t]}
        B ( \Delta A_u,\Delta X_u ).
\end{equation*}
This implies that
\begin{align*}
    &
    \lim_{n \to \infty} 
        \left\lVert 
            \sum_{I \in \pi_n} B(\delta A_t,\delta X_t)(I) e^1_{D_{\varepsilon}}(I) 
            - \sum_{0 < u \leq t} B( \Delta A_u,\Delta X_u )
        \right\rVert   \\
    & \qquad = 
        \left\lVert 
            \sum_{u \in D^{\varepsilon} \cap [0,t]} B ( \Delta A_u,\Delta X_u )
        \right\rVert
    \leq 
        \lVert B \rVert \sup_{u \in [0,t]} \lVert \Delta X_u \rVert \sum_{u \in D^{\varepsilon} \cap [0,t]} \lVert \Delta A_u \rVert.
\end{align*}

Next, we consider the second term on the right-hand side of \eqref{5.6c}.
Since $X$ has no jumps greater than $\varepsilon$ on $I$
whenever $e^2_{D_{ \varepsilon}}(I) = 1$,
we have the estimate
\begin{equation*}
    \left\lVert \delta X_t(I) \right\rVert_F e^2_{D_{ \varepsilon}}(I)
    = 
    \lVert \delta (X - J_{D_\varepsilon}(X))_t(I) \rVert
        e^2_{D_{ \varepsilon}}(I)
    \leq 
    O_t^{+}(X-J_{D_{\varepsilon}}(X);\pi_n) e^2_{D_{ \varepsilon}}(I).
\end{equation*}
Hence
\begin{align*}
    \left\lVert 
        \sum_{I \in \pi_n} B(\delta A_t,\delta X_t)(I) e^2_{D_{ \varepsilon}}(I) 
    \right\rVert_{G} 
    & \leq
    \lVert B \rVert 
        \sum_{I \in \pi_n} \lVert \delta A_t(I) \rVert_E \,
        \lVert \delta X_t(I) \rVert_F \,
        e^2_{D_{ \varepsilon}}(I)   \\
    & \leq \lVert B \rVert O_t^{+}(X-J_{D_{\varepsilon}}(X);\pi_n) V(A)_t.
\end{align*}

From the discussion above, we can deduce that
\begin{align*}
    & 
    \varlimsup_{n \to \infty} 
        \left\lVert 
            \sum_{I \in \pi_n} B(\delta A_t,\delta X_t )(I) 
            - \sum_{0 < u \leq t} B\left( \Delta A_u,\Delta X_u \right) 
        \right\rVert_{G}  \\
    & \qquad \leq 
    \lVert B \rVert \sup_{u \in [0,t]} \lVert \Delta X_u \rVert_F 
        \sum_{u \in D^{\varepsilon} \cap [0,t]} \lVert \Delta A_u \rVert_E
    + \lVert B \rVert V(A)_t \varlimsup_{n \to \infty} O_t^{+}(X-J_{D_{\varepsilon}}(X);\pi_n).
\end{align*}
Consequently, 
\begin{equation*}
    \varlimsup_{n \to \infty} 
        \left\lVert 
            \sum_{I \in \pi_n} B(\delta A_t,\delta X_t)(I)
            - \sum_{0 < u \leq t} B\left( \Delta A_u,\Delta X_u \right) 
        \right\rVert_{G}
    = 0,
\end{equation*}
which is the desired conclusion.
\end{proof}

Applying Proposition~\ref{5.6b} to the canonical 
bilinear maps $\otimes \colon E \times F \to E \widehat{\otimes}_{\alpha} F$
and $\otimes \colon F \times E \to F \widehat{\otimes}_{\alpha'} E$,
we get the following corollary.

\begin{corollary} \label{5.6d}
Let $A \in FV(\mathbb{R}_{\geq 0},E)$
and $X \in D(\mathbb{R}_{\geq 0},F)$.
If $(\pi_n)$ satisfies Condition (C) for $(A,X)$,
then it has tensor quadratic covariations
${^\alpha [A,X]}$ and ${^{\alpha'} [X,A]}$ given by
\begin{equation*}
    {^\alpha [A,X]_t} = \sum_{0 < s \leq t} \Delta A_s \otimes \Delta X_s, \qquad
    {^{\alpha'} [X,A]_t} = \sum_{0 < s \leq t} \Delta X_s \otimes \Delta A_s
\end{equation*}
for every reasonable crossnorms $\alpha$ and $\alpha'$ on $E \otimes F$ and $F \otimes E$,
respectively.
\end{corollary}

Using Corollaries~\ref{5.4eb} and \ref{5.6d}, we obtain the following.

\begin{corollary} \label{5.6e}
Let $(X,A) \in D(\mathbb{R}_{\geq 0},E) \times FV(\mathbb{R}_{\geq 0},E)$
and suppose that $(\pi_n)$ satisfies Condition (C) for $(X,A)$.
If $X \colon \mathbb{R}_{\geq 0} \to E$ has
$\alpha$-tensor quadratic variation along $(\pi_n)$,
then $X+A$ has $\alpha$-tensor quadratic variation given by
\begin{equation*}
    {^\alpha [X+A,X+A]} = {^\alpha [X,X]} + {^\alpha [X,A]} + {^\alpha [A,X]} + {^\alpha [A,A]}
\end{equation*} 
for every reasonable crossnorm $\alpha$ on $E \otimes E$.
\end{corollary}

By a discussion similar to the proof of Proposition~\ref{5.6b},
we see that a path of finite variation has scalar quadratic variation.

\begin{proposition} \label{5.6f}
Let $A$ be a c\`{a}dl\`{a}g path of finite variation in a Banach space $E$.
If $(\pi_n)$ satisfies Condition (C) for $A$,
then $A$ has the scalar quadratic covariation given by 
\begin{equation*}
    Q(A)_t = \sum_{0 < s \leq t} \lVert \Delta A_s \rVert^2.
\end{equation*}
\end{proposition}

In the preceding part of this paper,
we have used the summation
\begin{equation*}
    \sum_{\lrbrack r,s \rbrack \in \pi_n}
        B(X_{s \wedge t}-X_{r \wedge t},Y_{s \wedge t}-Y_{r \wedge t})
\end{equation*}
to define the quadratic covariation.
We can also consider a different summation
\begin{equation*}
    \sum_{\substack{\lrbrack r,s \rbrack \in \pi_n}} 1_{\lbrack 0,t \rlbrack}(r) B(X_{s}-X_{r},Y_{s}-Y_{r}),
\end{equation*}
which is a slight modification of
that used in the original paper by F\"{o}llmer~\cite{Foellmer_1981}.
Let us investigate the relation between
these two summations.

\begin{proposition} \label{5.4f}
Let $(X,Y) \in D(\mathbb{R}_{\geq 0}, E \times F)$
and $B \in \mathcal{L}^{(2)}(E,F;G)$.
Suppose that $(X,Y)_{\overline{\pi_n}(t)} \to (X,Y)_t$
holds for all $t \in \mathbb{R}_{\geq 0}$.
Then the following two conditions are equivalent, respectively.
\begin{enumerate}
    \item The path $(X,Y)$ has strong/weak $B$-quadratic covariation along $(\pi_n)$.
    \item There exists a c\`{a}dl\`{a}g path $V \in FV(\mathbb{R}_{\geq 0},G)$ such that
        \begin{enumerate}
            \item for all $t \in \mathbb{R}_{\geq 0}$
                \begin{equation*}
                    \sum_{\lrbrack r,s \rbrack \in \pi_n}
                        1_{\lbrack 0,t \rlbrack}(r) B(\delta X,\delta Y)(\lrbrack r,s \rbrack)
                    \xrightarrow[n \to \infty]{} 
                    V_t
                \end{equation*}
                in the norm/weak topology of $G$,
            \item for all $t \in \mathbb{R}_{\geq 0}$
                \begin{equation*}
                    \Delta V_t = B(\Delta X_t,\Delta Y_t).
                \end{equation*}
        \end{enumerate}
\end{enumerate}
If these conditions are satisfied, 
then $V$ coincides with the quadratic covariation $Q_B(X,Y)$.
\end{proposition}

\begin{proof}
We show the assertion about strong convergence.
For each $t > 0$
\begin{align*}
    & 
        \left\lVert
            \sum_{\lrbrack r,s \rbrack \in \pi_n} B( \delta X_t,\delta Y_t)(\lrbrack r,s \rbrack)
            - \sum_{\lrbrack r,s \rbrack \in \pi_n} 1_{\lbrack 0,t \rlbrack}(r) B(\delta X,\delta Y)(\lrbrack r,s \rbrack)
        \right\rVert  \\
    & \qquad =
        \left\lVert 
            B \left( \delta X_t,\delta Y_t \right)(\pi_n(t))
            - B\left( \delta X,\delta Y \right)(\pi_n(t))
        \right\rVert    \\
    & \qquad \leq
        \left\lVert
            B \left( \delta X_t-\delta X, \delta Y_t \right)(\pi_n(t))
        \right\rVert
        + \left\lVert
            B \left( \delta X,\delta Y_t-\delta Y \right)(\pi_n(t))
        \right\rVert \\
    & \qquad \leq 
        \lVert B \rVert
            \lVert X_{\overline{\pi_n}(t)} - X_t \rVert
            \lVert \delta Y_t (\pi_n(t)) \rVert
        + \lVert B \rVert
            \lVert \delta X(\pi_n(t)) \rVert
            \lVert Y_{\overline{\pi_n}(t)} - Y_t \rVert.
\end{align*}
Hence, by assumption, the convergences of these two sequences are equivalent and
their limits coincide.

In the weak case, we have a similar estimate for the pairing
$\langle z^\ast,\phantom{x} \rangle$, which shows the assertion about
the weak quadratic covariation.
\end{proof}

According to Proposition~\ref{5.4f},
we see that the two definitions of quadratic covariation are equivalent
provided that $(\pi_n)$ satisfies the assumption in the proposition.
The first definition,
which is given in Definition~\ref{2b}, is more intuitive.
The second one has some technical advantages because the path
$
    t \mapsto 
    \sum_{I \in \pi_n} 
        1_{\lbrack 0,t \rlbrack} B(\delta X,\delta Y)(I)$
is of finite variation
\footnote{
    Note that the path 
    $
        t \mapsto 
        \sum_{\pi_n} 1_{\lbrack 0,t \rlbrack} B(\delta X,\delta Y)(I)
    $
    is \emph{c\`{a}gl\`{a}d} but not c\`{a}dl\`{a}g.
}.

\begin{remark} \label{5.5c}
Following a discussion similar to that in Proposition~\ref{5.4f},
we can obtain an equivalent definition of scalar quadratic variation using the summation
$
    \sum_{I \in \pi_n} 
        \lVert 1_{\lbrack 0,t \rlbrack} \delta X (I) \rVert^2
$
if $(\pi_n)$ satisfies the same condition as Proposition~\ref{5.4f}.
\end{remark}

\section{Proof of Lemma \ref{2g}}

In this section, we prove Lemma~\ref{2g},
which is essentially used to show the main theorems of this paper.
To prove it, we prepare some additional lemmas.
Though Lemma~\ref{2g} includes both weak and strong convergence results,
we mainly focus on the proof of the weak case
\footnote{
    For the strong case, the reader can refer to an earlier version of
    this article at arXiv: \href{https://arxiv.org/abs/2104.08138v2}{2104.08138v2}.
}.

Throughout this section, let the symbols $E$, $E_1$, and $G$ denote Banach spaces
and $B \colon E \times E \to E_1$ denote a bounded bilinear map.

\begin{lemma} \label{5.7b}
Let $s > 0$ and let $(\pi_n)_{n \in \mathbb{N}}$ be a sequence of partitions.
Then $(\pi_n)$ approximates $1_{\lbrack 0,s \rlbrack}$ from the left if and only if
$\overline{\pi_n}(s) \to s$ as $n \to \infty$. 
\end{lemma}

\begin{proof}
First, assume that $(\pi_n)$ approximates $1_{\lbrack 0,s \rlbrack}$ from the left
and take an $\varepsilon > 0$ arbitrarily.
By definition, $1_{\lbrack 0,s \rlbrack}(s+\varepsilon-) = 0$.
Since $(\pi_n)$ approximates $1_{\lbrack 0,s \rlbrack}$ from the left
and $1_{\lbrack 0,s \rlbrack}$ takes only two values $0$ and $1$,
one sees that $1_{\lbrack 0,s \rlbrack}(\underline{\pi_n}(s+\varepsilon)) = 0$ holds
for a large enough $n$.
Therefore, $\underline{\pi_n}(s+\varepsilon) \geq s$ for large enough $n$.
This leads to 
$\overline{\pi_n}(s) \leq \underline{\pi_n}(s+\varepsilon) < s + \varepsilon$
for a sufficiently large $n$.
Hence $\overline{\pi_n}(s) \to s$ as $n \to \infty$.

Conversely, assume that $(\overline{\pi_n}(s))_{n \in \mathbb{N}}$
converges to $s$ as $n \to \infty$.
If $0 < t \leq s$, then the convergence 
$1_{\lbrack 0,s \rlbrack}(\underline{\pi_n}(t)) \to 1_{\lbrack 0,s \rlbrack}(t-) = 1$ is obvious.
Let $t > s$. We see from the assumption that $\overline{\pi_n}(s) < t$ for large enough $n$.
For such an $n$, we have $s \leq \overline{\pi_n}(s) \leq \underline{\pi_n}(t) < t$,
and hence $1_{\lbrack 0,s \rlbrack}(\underline{\pi_n}(t)) = 0 = 1_{\lbrack 0,s \rlbrack}(t-)$.
This shows that $\lim_{n \to \infty} 1_{\lbrack 0,s \rlbrack}(\underline{\pi_n}(t)) = 1_{\lbrack 0,s \rlbrack}(t-)$.
Hence $(\pi_n)$ is a left approximation sequence for $1_{\lbrack 0,s \rlbrack}$.
\end{proof}

\begin{lemma} \label{5.7c}
Suppose that $X \in D(\mathbb{R}_{\geq 0},E)$ has
weak $B$-quadratic variation along a sequence $(\pi_n)$
satisfying Condition (C) for $X$.
Let $r,s$ be two real numbers satisfying $0 \leq r < s$.
If $(\pi_n)$ approximates $1_{\lbrack r,s \rlbrack}$
from the left, then 
\begin{equation} \label{5.7cb}
    \lim_{n \to \infty} \sum_{I \in \pi_n} 
        1_{\lbrack r,s \rlbrack} B(\delta X_t,\delta X_t)(I)
    = 
    Q_B(X,X)_{s \wedge t} - Q_B(X,X)_{r \wedge t}
\end{equation}
holds for all $t \in \mathbb{R}_{\geq 0}$
in the weak topology.

If $X$ has strong $B$-quadratic variation, then the convergence of \eqref{5.7cb}
holds in the norm topology.
\end{lemma}

\begin{proof}
We show the case of weak convergence.
By considering the decomposition 
\begin{equation*}
    \sum_{I \in \pi_n} 1_{\lbrack r,s \rlbrack} B(\delta X_t,\delta X_t)(I)
    = 
    \sum_{I \in \pi_n} 1_{\lbrack 0,s \rlbrack} B(\delta X_t,\delta X_t)(I)
        - \sum_{I \in \pi_n} 1_{\lbrack 0,r \rlbrack} B(\delta X_t,\delta X_t)(I),
\end{equation*}
we can assume that $r=0$ without loss of generality.

If $t \leq s$, the equation
\begin{equation*}
    \sum_{I \in \pi_n} 1_{\lbrack 0,s \rlbrack} B(\delta X_t,\delta X_t)(I)
    =
    \sum_{I \in \pi_n} B(\delta X_t,\delta X_t)(I)
\end{equation*}
holds for all $n \in \mathbb{N}$.
Therefore,
\begin{align*}
    \lim_{n \to \infty}
        \left\langle 
            z^*,
            \sum_{I \in \pi_n} 1_{\lbrack 0,s \rlbrack} B(\delta X_t,\delta X_t)(I)
        \right\rangle
    =
    \langle z^*,Q_B(X,X)_{t} \rangle
    =
    \langle z^*,Q_B(X,X)_{t \wedge s} \rangle
\end{align*}
for all $z^* \in E_1^*$.

Next, assume $s < t$.
Then, by Lemma~\ref{5.7b}, $\overline{\pi_n}(s) \to s$ as $n \to \infty$.
Hence,
\begin{align*}
    \lim_{n \to \infty}
        \left\langle 
            z^*,
            \sum_{I \in \pi_n} 1_{\lbrack 0,s \rlbrack} B(\delta X_t,\delta X_t)(I)
        \right\rangle
    & =
    \lim_{n \to \infty}
        \left\langle 
            z^*,
            \sum_{I \in \pi_n} 1_{\lbrack 0,s \rlbrack} B(\delta X,\delta X)(I)
        \right\rangle   \\
    & =
    \langle z^*,Q_B(X,X)_s \rangle
    =
    \langle z^*,Q_B(X,X)_{s \wedge t} \rangle.
\end{align*}
Note that the second equality follows from the same argument
as the proof of Proposition~\ref{5.4f}.

In both cases, we have the desired convergence.

If $X$ has strong $B$-quadratic variation, we can directly show the norm convergence
of the sequence without taking the pairing $\langle z^\ast, \phantom{x} \rangle$.
\end{proof}

\begin{lemma} \label{5.7d}
Let $X$ be a c\`{a}dl\`{a}g path in $E$
with weak $B$-quadratic variation
along a sequence $(\pi_n)$ satisfying Condition (C) for $X$.
Suppose that 
$\xi \in D(\mathbb{R}_{\geq 0}, \mathcal{L}(E_1,G))$
has the representation
\begin{equation} \label{5.7e}
    \xi = \sum_{i \geq 1} 1_{\lbrack \tau_{i-1},\tau_{i} \rlbrack} a_{i},
\end{equation}
where $0 = \tau_0 < \tau_1 < \dots < \tau_i < \tau_{i+1} < \dots \to \infty$
and each $a_i$ is an element of $\mathcal{L}(E_1,G)$.
If $(\pi_n)$ approximates $\xi$ from the left,
then the Stieltjes integral of $\xi_{-}$ with respect
to $Q_B(X,X)$ is approximated as
\begin{equation} \label{5.7f}
    \lim_{n \to \infty} \sum_{I \in \pi_n} 
        \xi B(\delta X_t,\delta X_t)(I)
    = \int_{]0,t]} \xi_{s-}\mathrm{d}Q_B(X,X)_s.
\end{equation}
in the weak topology.

If $X$ has strong $B$-quadratic variation,
then \eqref{5.7f} holds in the norm topology of $G$.
\end{lemma}

\begin{proof}
We show the weak convergence case.
First, note that the Stieltjes integral on the right-hand side 
of \eqref{5.7f} has the representation
\begin{equation*}
    \int_{]0,t]} \xi_{s-}\mathrm{d}Q_B(X,X)_s
    = 
    \sum_{i \geq 1} a_i \{ Q_B(X,X)_{\tau_i \wedge t} - Q_B(X,X)_{\tau_{i-1} \wedge t} \}.
\end{equation*}
On the other hand,
the summation on the left-hand side of \eqref{5.7f} is calculated as
\begin{align*}
    \sum_{I \in \pi_n} \xi B(\delta X_t,\delta X_t) (I)
    = 
    \sum_{i \geq 1}  
        a_{i}
        \sum_{I \in \pi_n} 1_{\lbrack \tau_{i-1},\tau_{i} \rlbrack} 
                B(\delta X_t,\delta X_t)(I).
\end{align*}
Therefore, it suffices to show that 
\begin{equation*}
    \lim_{n \to \infty}
        \left\langle 
            z^*a_i,
            \sum_{I \in \pi_n} 1_{\lbrack \tau_{i-1},\tau_{i} \rlbrack} B(\delta X_t,\delta X_t)(I)
        \right\rangle
    =
    \langle z^*a_i, Q_B(X,X)_{\tau_i \wedge t} - Q_B(X,X)_{\tau_{i-1} \wedge t} \rangle
\end{equation*}
for all $z^* \in G^*$ and $i \geq 1$.
This follows directly from Lemma~\ref{5.7c}.

If $Q_B(X,X)$ is the strong $B$-quadratic variation,
the sequence of discrete sums converges in the norm topology
by the strong version of Lemma~\ref{5.7c}.
\end{proof}

\begin{lemma} \label{5.7g}
Let $V$ be a locally convex Hausdorff topological vector space
of which topology is generated by the family of 
seminorms $(\rho_{i})_{i \in I}$.
\begin{enumerate}
    \item Let $f \colon \mathbb{R}_{\geq 0} \to V$ be a c\`{a}dl\`{a}g path.
        Then for every $i \in I$ and $\varepsilon > 0$, 
        there is a right continuous step function $h$ such that
        $\rho_i(f(t) - h(t)) \leq \varepsilon$ for all $t \geq 0$.
    \item If a sequence of partitions $(\pi_n)$ approximates $f$ from the left,
        then the step function $h$ in (i) can be chosen so that 
        $(\pi_n)$ still approximates $h$ from the left.
\end{enumerate}
\end{lemma}

\begin{proof}
Fix $\varepsilon > 0$ and $i \in I$.
We define the oscillation of $f$
on $S \subset \mathbb{R}_{\geq 0}$ by
\begin{equation*}
    \omega_i(f;S) = \sup \{ \rho_i(f(t) - f(s)) \mid s,t \in S \}.
\end{equation*}

(i)
Let us construct a partition of $\mathbb{R}_{\geq 0}$ recursively.
First, let $t_0 = 0$.
Next, assume that there is a sequence $0 = t_0 < \dots < t_n$
such that $\omega_i(f; \lbrack t_k,t_{k+1} \rlbrack) \leq \varepsilon$ for $k \in \{ 0,\dots, n-1 \}$.

\emph{Case A}.
If $\omega_i(f;\lbrack t_n,\infty \rlbrack) = 0$, then we set $t_{n+1} = t_n + 1$.

\emph{Case B}.
If $\omega_i(f; \lbrack t_n,\infty \rlbrack) > \varepsilon$,
first let
\begin{equation*}
    t_{n+1}' 
    = 
    \inf \left\{ 
        t > t_n \, \middle\vert\,  \omega_i(f; \lbrack t_n, t \rlbrack) > \varepsilon 
    \right\},
\end{equation*}
and then define
\begin{equation} \label{5.7gb}
    t_{n+1}
    = 
    \begin{cases}
        \sup \{ 
            t > t_{n+1}' \mid  \omega_i(f; \lbrack t_{n+1}', t \rlbrack) = 0
        \} & \text{if} \ \omega_i(f; \lbrack t_{n+1}', \infty \rlbrack) > 0  \\
        t_{n+1}' & \text{otherwise}.
    \end{cases}
\end{equation}
Note that $f$ is not constant on any interval of the form
$\lbrack t_{n+1}, t \rlbrack$ ($t > t_{n+1}$)
when $\omega_i(f; \lbrack t_{n+1}', \infty \rlbrack) > 0$.

\emph{Case C}.
If $0 < \omega_i(f; \lbrack t_n,\infty \rlbrack) \leq \varepsilon$,
then let $t_{n+1}' = t_n + 1$ and define $t_{n+1}$ by the formula~\eqref{5.7gb}.

The sequence defined above satisfies $t_n \to \infty$.
Indeed, if $\omega_i( f;[t_n,\infty \rlbrack) \leq \varepsilon$ for some $n \in \mathbb{N}$,
we have $t_{n+k} \geq t_n + k$ for all $k \in \mathbb{N}$, by definition.
Now assume that $\omega_i(f;[t_n,\infty \rlbrack) > \varepsilon$ for all $n \in \mathbb{N}$
and $t_n \uparrow t^* < \infty$.
In this case, $\omega_i(f;[t^*-\delta,t^* \rlbrack) \geq \varepsilon$ holds for arbitrary small $\delta$.
This contradicts the existence of the left limit at $t^*$.

Let $P \coloneqq (t_n)_{n \in \mathbb{N}}$ and define a function $f_{P}$ by the formula
\begin{equation*}
    f_{P} = \sum_{n \geq 0} 1_{\lbrack t_n,t_{n+1} \rlbrack} f(t_{n+1}-).
\end{equation*}
For $t \in \lbrack t_n,t_{n+1} \rlbrack$, we see that
\begin{equation*}
    \rho_i(f(t) - f_P(t)) 
    = \rho_i(f(t) - f(t_{n+1}-))
    \leq \omega_i(f;\lbrack t_n,t_{n+1} \rlbrack) \leq \varepsilon.
\end{equation*}
Hence $f_P$ is a right continuous step function satisfying the desired condition.

(ii)
We shall show that the function $f_P$ defined above is approximated from left 
by the left-approximation sequence $(\pi_k)$ for $f$.
Note that, by the definition of $f_P$, 
\begin{equation*}
    f_P(s-) = \sum_{n \geq 0} 1_{\lrbrack t_n,t_{n+1} \rbrack}(s) f(t_{n+1}-)
\end{equation*}
holds for all $s > 0$.

Now fix $s > 0$ arbitrarily and choose a unique $n_0 \in \mathbb{N}$
such that $s \in \lrbrack t_{n_0},t_{n_0+1} \rbrack$.
If $f$ is discontinuous at $t_{n_0}$, then we see that
$\overline{\pi_k}(t_{n_0}) \to t_{n_0}$ by the same argument
as the proof of Lemma~\ref{5.7b}.
In this case, $t_{n_0} \leq \overline{\pi_k}(t_{n_0}) \leq \underline{\pi_k}(s) < s$ for large enough
$k$ and therefore we have
\begin{equation*}
    \lim_{k \to \infty} f_P(\underline{\pi_k}(s)) = f(t_{n_0+1}-) = f_P(s-).
\end{equation*}

Next, assume that $f$ is continuous at $t_{n_0}$ and $f$ is not constant
on any interval of the form $\lbrack t_{n_0},t \rlbrack$ for $t > t_{n_0}$.
Take $t'$ such that $t_{n_0} < t' < s$ and $f(t'-) \neq f(t_{n_0}) = f(t_{n_0}-)$.
If $\underline{\pi_k}(t') \leq t_{n_0}$ for infinitely many $k$,
then we can take a subsequence such that
$\underline{\pi_{k_l}}(t') \leq t_{n_0}$ for all $l$.
By definition, $\underline{\pi_{k_l}}(t') = \underline{\pi_{k_l}}(t_{n_0})$
and hence
\begin{equation*}
    f(t'-)
    = \lim_{l \to \infty} f(\underline{\pi_{k_l}}(t'))
    = \lim_{l \to \infty} f(\underline{\pi_{k_l}}(t_{n_0}))
    = f(t_{n_0}-),
\end{equation*}
which contradicts the assumption on $t'$.
This shows that
$t_{n_0} < \underline{\pi_k}(t') \leq \underline{\pi_k}(s) < s$
holds for large enough $k$ and consequently
\begin{equation*}
    \lim_{k \to \infty} f_P(\underline{\pi_k}(s)) = f(t_{n_0+1}-) = f_P(s-).
\end{equation*}

Finally, assume that $f$ is continuous at $t_{n_0}$
and $f$ is constant on $\lbrack t_{n_0},t \rlbrack$ for some $t > t_{n_0}$.
In this case, we have $\omega_i(f;\lbrack t_{n_0},\infty \rlbrack) = 0$
by the definition of $P = (t_n)$.
Also note that $t_{n_0} \leq \underline{\pi_k}(s) < s$ or
$\underline{\pi_k}(t_{n_0}) = \underline{\pi_k}(s)$
holds for each $k$. In both cases,
\begin{equation*}
    \lim_{k \to \infty} f_P(\underline{\pi_k}(s))
    = f(t_{n_0}-)
    = f(t_{n_0+1}-)
    = f_P(s-).
\end{equation*}
The second equality follows from the continuity of $f$ at $t_{n_0}$
and the property that $\omega_i(f;\lbrack t_{n_0},\infty \rlbrack) = 0$.
This completes the proof.
\end{proof}

Finally, we start dealing with the proof of Lemma~\ref{2g}.

\begin{proof}[Proof of Lemma \ref{2g}]
First, fix $t > 0$ and $\varepsilon > 0$ arbitrarily
and choose a compact set $K \subset E_1$ satisfying
\begin{equation*}
    B(\delta X_t,\delta X_t)(\mathcal{I} \cap [0,t]) \cup \delta Q_B(X,X)_t(\mathcal{I} \cap [0,t]) \subset K.
\end{equation*}
By Lemma~\ref{5.7g}, we can find an 
$\mathcal{L}_{\mathrm{c}}(E_1,G)$-valued right continuous step function
\begin{equation*}
    h = \sum_{i \geq 1} 1_{\lbrack \tau_{i-1},\tau_i \rlbrack} a_i
\end{equation*}
so that $\rho_K(h(s) - \xi(s)) \leq \varepsilon$ holds for all $s \in [0,t]$
and $(\pi_n)$ approximates $h$ from the left.
Then, for all $z^* \in G^*$,
\begin{align} \label{5.7h}
    &
    \left\lvert 
        \left\langle
            z^*,
            \sum_{I \in \pi_n}
                \xi B(\delta X_t,\delta X_t )(I)
        \right\rangle
        - 
        \left\langle 
            z^*,
            \int_{\lrbrack 0,t \rbrack} \xi_{u-} \mathrm{d}Q_B(X,X)_u
        \right\rangle
    \right\rvert   \notag \\
    & \qquad \leq 
    \left\lvert
        \sum_{I \in \pi_n}
            \left\langle 
                z^*\xi, B(\delta X_t,\delta X_t )
            \right\rangle(I)
        - \sum_{I \in \pi_n} 
            \left\langle
                z^* h, B(\delta X_t,\delta X_t )
            \right\rangle(I)
    \right\rvert      \notag \\ 
    & \qquad\quad +
    \left\lvert
        \left\langle 
            z^*,
            \sum_{\lrbrack r,s \rbrack \in \pi_n}
                h B(\delta X_t,\delta X_t )(I)
        \right\rangle
        - 
        \left\langle 
            z^*,
            \int_{\lrbrack 0,t \rbrack} h(u-) \mathrm{d}Q_B(X,X)_u
        \right\rangle 
    \right\rvert     \notag  \\
    & \qquad\quad +
    \left\lvert
        \int_{\lrbrack 0,t \rbrack}
            \left\langle
                z^* h(u-), \mathrm{d}Q_B(X,X)_u
            \right\rangle 
        - \int_{\lrbrack 0,t \rbrack}
            \left\langle
                z^* \xi(u-), \mathrm{d}Q_B(X,X)_u 
            \right\rangle
    \right\rvert.
\end{align}

We will observe the behaviour of each part of the right-hand side.
We can deduce from Lemma~\ref{5.7d} that
the second term converges to $0$ as $n \to \infty$.
By the choice of $h$, we find that
\begin{align*}
    \left\lvert
    \sum_{I \in \pi_n}
        \left\langle 
            z^* \xi, B(\delta X_t,\delta X_t)
        \right\rangle(I)
    - \sum_{I \in \pi_n} 
        \left\langle
            z^* h, B(\delta X_t,\delta X_t)
        \right\rangle(I)
    \right\rvert
    & \leq
    \varepsilon \lVert z^* \rVert \lVert B \rVert
    \sum_{I \in \pi_n}
        \left\lVert 
            \delta X_t(I)
        \right\rVert_E^2.
\end{align*}
Therefore,
\begin{equation*}
    \varlimsup_{n \to \infty} 
        \left\lvert
        \sum_{I \in \pi_n}
            \left\langle 
                z^*\xi, B(\delta X_t,\delta X_t)
            \right\rangle(I)
        - \sum_{I \in \pi_n} 
            \left\langle
                z^*h, B(\delta X_t,\delta X_t)
            \right\rangle(I)
        \right\rvert
    \leq \varepsilon \lVert z^* \rVert \lVert B \rVert V^2(X;\Pi)_t.
\end{equation*}
On the other hand, we have
\begin{align*}
    \left\lvert
        \int_{\lrbrack 0,t \rbrack}
            \left\langle
                z^* ( h(u-)-\xi(u-) ), \mathrm{d}Q_B(X,X)_u 
            \right\rangle
    \right\rvert
& \leq 
    \lVert z^* \rVert
        \int_{\lrbrack 0,t \rbrack} \rho_K (h(u-)-\xi(u-)) \mathrm{d}V(Q_B(X,X))_t  \\
& \leq
    \varepsilon \lVert z^* \rVert V(Q_B(X,X))_t
\end{align*}
by Proposition~\ref{a3f}.
Consequently,
\begin{align*}
    \varlimsup_{n \to \infty}
    \left\lvert 
        \left\langle 
            z^*,
            \sum_{I \in \pi_n} \xi B(\delta X_t, \delta X_t)(I)
            - \int_{\lrbrack 0,t \rbrack} \xi_{u-} \mathrm{d}Q_B(X,X)_u
        \right\rangle
    \right\rvert
    \leq
    \varepsilon \lVert z^* \rVert \{ V^2(X;\Pi)_t + V(Q_B(X,X))_t \}.
\end{align*}
Because $\varepsilon$ is chosen arbitrarily,
we get the desired conclusion.

If $Q_B(X,X)$ is the strong quadratic variation, we replace \eqref{5.7h}
with a similar norm inequality.
In this case, we see that the corresponding second term converges to $0$
by the strong version of Lemma~\ref{5.7d}.
The remaining terms are estimated in almost the same way as above.
\end{proof}

\section{The It{\^o} formula}

This section is devoted to showing the It{\^o} formula
within our framework of the It\^{o}--F\"{o}llmer calculus in Banach spaces.
Let us begin by defining It\^{o}--F\"{o}llmer integrals.

\begin{definition} \label{5.9b}
Let $E$ be a locally convex space, and let $F$ and $G$ be Banach spaces.
Consider c\`{a}dl\`{a}g paths $H$ and $X$ in $E$ and $F$, respectively,
and a continuous bilinear map $B \colon E \times F \to G$.
Suppose that a sequence of partitions $(\pi_n)$ approximates
$H$ from the left. We call the limit
\begin{equation*}
    \int_0^t B(H_{s-}, \mathrm{d}X_{s})
    = 
    \int_{\lrbrack 0,t \rbrack} B(H_{s-}, \mathrm{d}X_{s})
    \coloneqq 
    \lim_{n \to \infty} \sum_{\lrbrack r,s \rbrack \in \pi_n} 
        B\left( H_r, \delta X_t(\lrbrack r,s \rbrack) \right) \in G
\end{equation*}
the
\emph{%
    (strong) It\^{o}--F\"{o}llmer integral of $H$ with respect to $X$ along $(\pi_n)$%
}
if it exists.
If this convergence holds in the weak topology, 
we call the limit the \emph{weak It\^{o}--F\"{o}llmer integral}.

Similarly, the strong and the weak It\^{o}--F\"{o}llmer integral for a $B' \in \mathcal{L}^{(2)}(F,E;G)$
are defined as the limit 
\begin{equation*}
    \int_{\lrbrack 0,t \rbrack} B'(\mathrm{d}X_s, H_{s-})
    = 
    \lim_{n \to \infty} \sum_{\lrbrack r,s \rbrack \in \pi_n} 
        B'\left( \delta X_t(\lrbrack r,s \rbrack),H_{r} \right) \in G
\end{equation*}
with the corresponding topology.
\end{definition}

If $B$ is the canonical bilinear map
$\otimes \colon E \times F \to E \widehat{\otimes}_{\alpha} F$,
we write
\begin{equation*}
    \int_{\lrbrack 0,t \rbrack} B(H_{s-},\mathrm{d}X_s)
    = 
    \int_{\lrbrack 0,t \rbrack} H_{s-} \otimes \mathrm{d}X_s.
\end{equation*}

\begin{remark} \label{5.9c}
The It\^{o}--F\"{o}llmer integral of Definition~\ref{5.9b}
inherits the bilinear property from
$B \in \mathcal{L}^{(2)}(E,F;G)$ in the following sense. 
\begin{enumerate}
    \item Suppose that the following two It\^{o}--F\"{o}llmer integrals exist:
        \begin{equation*}
            \int_{\lrbrack 0,t \rbrack} B(H_{s-}, \mathrm{d}X_{s}), \qquad 
            \int_{\lrbrack 0,t \rbrack} B(K_{s-}, \mathrm{d}X_{s}).
        \end{equation*}
        Then, for every $\alpha$, $\beta \in \mathbb{R}$,
        the It\^{o}--F\"{o}llmer integral of $\alpha H + \beta K$
        with respect to $X$ exists and satisfies
        \begin{equation*}
            \int_{\lrbrack 0,t \rbrack} B(H_{s-} + K_{s-}, \mathrm{d}X_{s})
            =
            \alpha \int_{\lrbrack 0,t \rbrack} B(H_{s-}, \mathrm{d}X_{s}) 
            + \beta \int_{\lrbrack 0,t \rbrack} B(K_{s-}, \mathrm{d}X_{s}).
        \end{equation*}
    \item Suppose that the following two It\^{o}--F\"{o}llmer integrals exist:
        \begin{equation*}
            \int_{\lrbrack 0,t \rbrack} B(H_{s-}, \mathrm{d}X_{s}), \qquad 
            \int_{\lrbrack 0,t \rbrack} B(H_{s-}, \mathrm{d}Y_{s}).
        \end{equation*}
        Then, for every $\alpha$, $\beta \in \mathbb{R}$,
        the It\^{o}--F\"{o}llmer integral of $H$ with respect to 
        $\alpha X + \beta Y$ exists and satisfies
        \begin{equation*}
            \int_{\lrbrack 0,t \rbrack} B(H_{s-}, \mathrm{d}(\alpha X+ \beta Y)_{s})
            =
            \alpha \int_{\lrbrack 0,t \rbrack} B(H_{s-}, \mathrm{d}X_{s}) 
            + \beta \int_{\lrbrack 0,t \rbrack} B(H_{s-}, \mathrm{d}Y_{s}).
        \end{equation*}
\end{enumerate}
\end{remark}

First, we consider the case where
the integrator is a path of finite variation.
From the dominated convergence theorem,
we can easily deduce the following proposition.

\begin{proposition} \label{5.9d}
Let $H \in D(\mathbb{R}_{\geq 0},E)$, $A \in FV(\mathbb{R}_{\geq 0},F)$,
and $B \in \mathcal{L}^{(2)}(E,F;G)$.
If a sequence of partitions $(\pi_n)$ approximates $H$ from the left,
we have
\begin{equation*}
    (\mathrm{IF}) \int_{\lrbrack 0,t \rbrack} B(H_{s-},\mathrm{d}A_s)
    =
    (\mathrm{S}) \int_{\lrbrack 0,t \rbrack} B(H_{s-},\mathrm{d}A_s).
\end{equation*}
Here, the integral of the left-hand side is the It\^{o}--F\"{o}llmer integral
by Definition~\ref{5.9b},
and that of the right-hand side is the usual Stieltjes integral.
\end{proposition}

Now we start to prove our main theorem.

\begin{proof}[Proof of Theorem \ref{2e}]
We show the weak convergence of the It\^{o}--F\"{o}llmer integral.

First, fix $t>0$ arbitrarily and choose compact convex sets
$K_1 \subset F$, $K_2 \subset E$,
and $K_3 \subset G$ such that
\begin{equation*}
    A([0,t]) \subset K_1, \quad 
    X([0,t]) \subset K_2, \quad 
    B(\delta X,\delta X)(\mathcal{I} \cap [0,t]) \subset K_3.
\end{equation*}

\emph{Step 1: Convergence of the summation in the formula~\eqref{2c}.}
In this step, we confirm that the summation of jump terms converges absolutely.
This is proved by the following estimate,
which follows from Taylor's formula (Proposition~\ref{a1c}):
\begin{align*}
    %1
    &
    \sum_{0 < s \leq t}
    \lVert
        f(A_s,X_s) - f(A_{s-},X_{s-})
        - D_x f(A_{s-},X_{s-})\Delta X_s
    \rVert  \\
    %2
    & \qquad \leq
    \sum_{0 < s \leq t}
    \lVert
        f(A_{s-},X_s) - f(A_{s-},X_{s-})
        - D_x f(A_{s-},X_{s-})\Delta X_s
    \rVert  + 
    \sum_{0 < s \leq t} \lVert f(A_s,X_s) - f(A_{s-},X_s) \rVert  \\
    & \qquad \leq 
    \sup_{(a,x) \in K_1 \times K_2} \lVert D_B^2 f(a,x) \rVert \lVert B \rVert \sum_{0 < s \leq t} \lVert \Delta X_s \rVert^2
    +
    \sup_{(a,x) \in K_1 \times K_2} \lVert D_a f(a,x) \rVert \sum_{0 < s \leq t} \lVert \Delta A_s \rVert  \\
    & \qquad < \infty.
\end{align*}
Note that the uniform boundedness principle
combined with the strong continuity
shows that $\sup_{(a,x)} \left\lVert D_B^2 f(a,x) \right\rVert$
and $\sup_{(a,x)} \lVert D_a f(a,x) \rVert$
are finite.

Moreover, we see that 
\begin{align*}
    \sum_{0 < s \leq t}
        \lVert
            D_x^2 f(A_{s-},X_{s-})\Delta X_s^{\otimes 2}
        \rVert
    & = 
    \sum_{0 < s \leq t}
        \lVert
            D_B^2 f(A_{s-},X_{s-})B(\Delta X_s^{\otimes 2}) 
        \rVert   \\
    & \leq
    \sup_{(a,x) \in K_1 \times K_2}
        \left\lVert D_B^2 f(a,x) \right\rVert \lVert B \rVert
        \sum_{0 < s \leq t} \left\lVert \Delta X_s \right\rVert^2
    < \infty
\end{align*}
and
\begin{equation*}
    \sum_{0 < s \leq t}
        \left\lvert
            D_a f(A_{s-},X_{s-})\Delta A_s
        \right\rvert
    \leq
    \sup_{(a,x) \in K_1 \times K_2} \lVert D_a f(a,x) \rVert \sum_{0 < s \leq t} \lVert \Delta A_s \rVert
    <
    \infty.
\end{equation*}
This shows that equation~\eqref{2f}
is equivalent to the following equation:
\begin{align} \label{5.9e}
    & 
    f(A_t,X_t) - f(A_0,X_0) 
    -
    \int_{\lrbrack 0,t \rbrack}
        D_x f (A_{s-},X_{s-})\, \mathrm{d}X_s   \notag \\
    & \qquad = 
    \int_{\lrbrack 0,t \rbrack} D_a f (A_{s-},X_{s-})\, \mathrm{d}A_s 
    +
    \frac{1}{2} \int_{\lrbrack 0,t \rbrack} 
        D_B^2 f (A_{s-},X_{s-})\, \mathrm{d}Q_B(X,X)_s     \notag \\
    & \qquad\quad +
    \sum_{0 < s \leq t} 
        \left\{ \Delta f (A_s,X_s) - D_x f (A_{s-},X_{s-})\Delta X_s    \right\}  \notag \\
    & \qquad\quad - 
    \sum_{ 0 < s \leq t} 
        D_a f (A_{s-},X_{s-})\Delta A_s
        - 
        \frac{1}{2} \sum_{0 < s \leq t} 
            D_B^2 f (A_{s-},X_{s-})B(\Delta X_s,\Delta X_s).
\end{align}
We will therefore prove \eqref{5.9e} instead of \eqref{2f}.

\emph{Step 2: The Taylor expansion.}
Let $I = \lrbrack u,v \rbrack \in \pi_n$ and 
consider the first-order Taylor expansion with respect to the variable $a$
on $I \cap [0,t] = \lrbrack u \wedge t, v \wedge t \rbrack$.
Then we have 
\begin{equation} \label{5.9f}
    f(A_{v \wedge t},X_{v \wedge t}) - f(A_{u \wedge t},X_{v \wedge t})
    =
    D_a f(A, X) \delta A_t(I)
    +
    r_t(I)
\end{equation}
where
\begin{equation*}
    r_t(I)
    =
    \int_{[0,1]}
        \left\{
            D_a f \left( A_{u \wedge t} + \theta \delta A_t(I), X_{v \wedge t} \right)
            -
            D_a f\left( A_{u \wedge t},X_{u \wedge t} \right)
        \right\} \delta A_t(I) \, 
    \mathrm{d}\theta.
\end{equation*}
Here, recall that the notation `$D_a f(A, X) \delta A_t$'
was introduced in the second paragraph of Section 3.

Next, we consider the second-order Taylor expansion
\begin{align} \label{5.9g}
        f(A_{u \wedge t},X_{v \wedge t}) - f(A_{u \wedge t},X_{u \wedge t}) 
    =
        D_x f(A,X) \delta X_t(I)
        + 
        \frac{1}{2} D_x^2 f(A,X) (\delta X_t)^{\otimes 2}(I)
        +
        R_t(I),
\end{align}
with $R_t(I)$ given by
\begin{align*}
    R_t(I)
    & =
    \frac{1}{2}
    \int_{[0,1]}
        (1-\theta)
        \left\{ 
            D_x^2 f(A_{u \wedge t},X_{u \wedge t} + \theta \delta X_t(I))
            -
            D_x^2 f(A_{u \wedge t},X_{u \wedge t})
        \right\} \delta X_t \otimes \delta X_t(I) \, 
    \mathrm{d}\theta   \\
    & =
    \frac{1}{2}
    \int_{[0,1]}
        (1-\theta)
        \left\{ 
            D_B^2 f(A_{u \wedge t},X_{u \wedge t} + \theta \delta X_t(I))
            -
            D_B^2 f(A_{u \wedge t},X_{u \wedge t})
        \right\} B(\delta X_t,\delta X_t)(I) \, 
    \mathrm{d}\theta.
\end{align*}
By the definition of $D_B^2 f$, we can rewrite \eqref{5.9g} as 
\begin{align} \label{5.9h}
        f(A_{u \wedge t},X_{v \wedge t}) - f(A_{u \wedge t},X_{u \wedge t})
    =
        D_x f(A,X)\delta X_t(I)
        + 
        \frac{1}{2} D_B^2 f(A,X)B(\delta X_t^{\otimes 2}) (I)
        +
        R_t(I),
\end{align}

Combining equations~\eqref{5.9g} and \eqref{5.9h}, we obtain
\begin{equation} \label{5.9hb}
    \delta f(A,X)_t
    =
    D_a f(A,X) \delta A_t
    + D_x f(A,X) \delta X_t 
    + \frac{1}{2} D_B f(A,X) B(\delta X_t^{\otimes 2})
    + r_t + R_t.
\end{equation}
Now take $\varepsilon > 0$ arbitrarily.
For simplicity, write $D = D(A,X)$,
$D_{\varepsilon} = D_{\varepsilon}(A,X)$,
and $D^{\varepsilon} = D^{\varepsilon}(A,X)$.
Using the notations $e^1_{D_{\varepsilon}}$ and $e^2_{D_{\varepsilon}}$
introduced in Section 5,
we can derive from \eqref{5.9hb} the equation
\begin{align*}
        \delta f(A,X)_t - \delta f(A,X)_t e^1_{D_{\varepsilon}} 
    & =
        D_x f(A,X) \delta X_t - D_x f(A,X) \delta X_t e^1_{D_{\varepsilon}}   \\ 
    & \quad +
        D_a f(A,X) \delta A_t - D_a f(A,X) \delta A_t e^1_{D_{\varepsilon}} + r_t e^2_{D_{\varepsilon}}   \\
    & \quad +
        \frac{1}{2} D_B^2 f(A,X) B(\delta X_t^{\otimes 2})
        - \frac{1}{2} D_B^2 f(A,X) B(\delta X_t^{\otimes 2}) e^1_{D_{\varepsilon}} 
        + R_t e^2_{D_{\varepsilon}}. 
\end{align*}
Moreover, by summing up this equality along $\pi_n$, we see that
\begin{align} \label{5.9i}
    &
    f(A_t,X_t)- f(A_0,X_0)
    -
    \sum_{I \in \pi_n} D_x f(A,X) \delta X_t(I)   \notag \\
%% 1
    & \qquad =
    \sum_{I \in \pi_n}
        \delta f(A,X)_t(I) e^1_{D_{\varepsilon}}(I)
%% 2
    -
    \sum_{I \in \pi_n} D_x f(A,X) \delta X_t(I) e^1_{D_{\varepsilon}}(I)  \notag \\
%% 3
    & \qquad\quad +
    \sum_{I \in \pi_n}
        D_a f(A,X) \delta A_t(I) 
%% 4
    -
    \sum_{I \in \pi_n}
        D_a f(A,X) \delta A_t(I) e^1_{D_\varepsilon}(I)  \notag \\
%% 5
    & \qquad\quad +
    \sum_{I \in \pi_n}
        \frac{1}{2} D_B^2 f(A,X) B(\delta X_t^{\otimes 2}) (I)
%% 6
    -
    \sum_{I \in \pi_n}
        \frac{1}{2} D_B^2 f(A,X) B(\delta X_t^{\otimes 2}) (I)
        e^1_{D_\varepsilon}(I) \notag \\
%% 7
    & \qquad\quad +
    \sum_{I \in \pi_n}
        r_t(I) e^2_{D_{\varepsilon}}(I)
%% 8
    +
    \sum_{I \in \pi_n}
        R_t(I) e^2_{D_{\varepsilon}}(I)  \notag \\
    & \qquad \eqqcolon
    I_1^{(n)}(t) - I_2^{(n)}(t) + I_3^{(n)}(t) - I_4^{(n)}(t)
    + I_5^{(n)}(t) - I_6^{(n)}(t) + I_7^{(n)}(t) + I_8^{(n)}(t).
\end{align}

\emph{Step 3: Behaviour of $I_1^{(n)}(t),\dots , I_8^{(n)}(t)$ of \eqref{5.9i}.}
Since $(\pi_n)$ satisfies Condition (C) for $(A,X)$, we can easily verify that
\begin{align*}
%%%
    \lim_{n \to \infty} I_4^{(n)}(t)
    & = 
        \sum_{s \in D_{\varepsilon} \cap [0,t]} 
            D_a f(A_{s-},X_{s-})\Delta A_s    \\
%%%
    \lim_{n \to \infty} I_2^{(n)}(t)
    & =
        \sum_{s \in D_\varepsilon \cap [0,t]} 
            D_x f(A_{s-},X_{s-}) \Delta X_s    \\
%%%
    \lim_{n \to \infty} I_6^{(n)}(t)
    & = 
        \frac{1}{2} \sum_{s \in D_\varepsilon \cap [0,t]} 
            D_x^2 f(A_{s-},X_{s-}) (\Delta X_s)^{\otimes 2}
    = 
        \frac{1}{2} \sum_{s \in D_\varepsilon \cap [0,t]} 
            D_B^2 f(A_{s-},X_{s-})B(\Delta X_s^{\otimes 2}).
\end{align*}
If $s \in D$ and $s \leq t$, we can deduce from Proposition~\ref{a1f} and the dominated convergence theorem that 
\begin{align*}
    \delta f(A,X)_t(\pi_n(s))
    & =
    \int_0^1 \{ Df[ (A,X)_{\underline{\pi_n}(s)} + \theta \delta (A,X)_t(\pi_n(s)) ] \} 
        \{ \delta (A,X)_t(\pi_n(s)) \} \,\mathrm{d} \theta   \\
    & \xrightarrow[n \to \infty]{}
    \int_0^1 Df[ (A,X)_{s-} + \theta \Delta (A,X)_s] \Delta (A,X)_s \,\mathrm{d} \theta  
    = \Delta f(A,X)_s
\end{align*}
Hence, 
\begin{equation*}
    \lim_{n \to \infty} I_1^{(n)}(t)
    = 
    \sum_{s \in D_{\varepsilon} \cap [0,t]} \Delta f (A_s,X_s).
\end{equation*}

By Lemma~\ref{2g}, we have
\begin{equation*}
    \lim_{n \to \infty} I_5^{(n)}(t)
    = 
    \frac{1}{2} \int_{\lrbrack 0,t \rbrack}  
        D_B^2 f(A_{s-},X_{s-})\mathrm{d}Q_B(X,X)_s
\end{equation*}
in the weak topology.

The dominated convergence theorem (Proposition~\ref{a3g}) gives
\begin{equation*}
    \lim_{n \to \infty} I_3^{(n)}(t)
    =
    \int_{\lrbrack 0,t \rbrack} D_a f(A_{s-},X_{s-}) \mathrm{d}A_s.
\end{equation*}

It remains to estimate the residual terms.
If $\lrbrack u,v \rbrack \cap D_{\varepsilon} = \emptyset$, then 
\begin{gather*}
    \omega(X, \lrbrack u,v \rbrack \cap [0,t])
        \leq 
            O^{+}_t(X-J(D_{\varepsilon}(X);X);\pi_n),   \\
    \omega(A, \lrbrack u,v \rbrack \cap [0,t])
        \leq 
            O^{+}_t(A-J(D_{\varepsilon}(A);A);\pi_n).
\end{gather*}
Now we write, for convenience,
\begin{gather*}
    \alpha(\varepsilon,n) = O^{+}_t(X-J(D_{\varepsilon}(X);X);\pi_n),   \\
    \beta(\varepsilon,n) = O^{+}_t(A-J(D_{\varepsilon}(A);A);\pi_n).
\end{gather*}
By the assumption that $(\pi_n)$ satisfies Condition (C) for $(A,X)$
--- and hence so does each of $A$ and $X$ --- we see that
\begin{equation*}
    \varlimsup_{\varepsilon \downarrow\downarrow 0} \varlimsup_{n \to \infty} \alpha(\varepsilon,n) = 0,  \qquad
    \varlimsup_{\varepsilon \downarrow\downarrow 0} \varlimsup_{n \to \infty} \beta(\varepsilon,n) = 0.
\end{equation*}
By definition,
\begin{align*}
    I_7^{(n)}(t)
    & \leq 
    \sum_{\lrbrack u,v \rbrack \in \pi_n}
        e^{2}_{D_{\varepsilon}}(\lrbrack u,v \rbrack)  \\
    & \quad \times 
        \int_{[0,1]}
            \lVert
                \left\{ D_a f \left( A_{u \wedge t} + \theta \delta A_t(\lrbrack u,v \rbrack), X_{v \wedge t} \right)
                -
                D_a f\left( A_{u \wedge t},X_{u \wedge t} \right) \right\} \delta A_t(\lrbrack u,v \rbrack)
            \rVert \,
        \mathrm{d}\theta   \\
    & \leq 
    V(A)_t
    \sup_{\substack{z,w \in K_1 \times K_2  \\ \lvert z-w \rvert \leq \alpha(\varepsilon,n) + \beta(\varepsilon,n)}}
        \rho_{K_1-K_1}(D_a f(z) - D_a f(w)).
\end{align*}
Similarly, we have 
\begin{align*}
    I_8^{(n)}(t)
    & \leq 
    \frac{1}{2}
    \sum_{\lrbrack u,v \rbrack \in \pi_n} e^{2}_{D_{\varepsilon}}(\lrbrack u,v \rbrack)  \\
    & \quad \times 
    \int_{[0,1]}
        (1-\theta)
        \bigl\lVert 
            \bigl\{ D_B^2 f(A_{u \wedge t},X_{u \wedge t} + \theta \delta X_t(\lrbrack u,v \rbrack))
            -
            D_B^2 f(A_{u \wedge t},X_{u \wedge t}) \bigr\} 
            B(\delta X_t^{\otimes 2})(\lrbrack u,v \rbrack)
        \bigr\rVert
    \mathrm{d}\theta   \\
    & \leq 
    \sup_{\substack{z,w \in K_1 \times K_2  \\ \lvert z-w \rvert \leq \alpha(\varepsilon,n)}}
        \rho_{K_3}(D_B^2 f(z) - D_B^2 f(w))
        \lVert B \rVert 
        \sum_{I \in \pi_n}\lVert \delta X_t(I) \rVert_{E}^2.
\end{align*}

Consequently, for every $z^* \in G^*$ with $\lVert z^* \rVert = 1$,
\begin{align} \label{5.9j}
%1
    &
    \varlimsup_{n \to \infty}
        \left\lvert
            \langle z^*,(\text{RHS of \eqref{5.9e}}) - (\text{RHS of \eqref{5.9i}}) \rangle
        \right\rvert    \notag \\
%2
    & \qquad \leq
    \left\lVert 
        \sum_{s \in D^{\varepsilon} \cap [0,t]}
            \left\{ 
                \Delta f(A_s,X_s) - \langle D_x f(A_{s-},X_{s-}), \Delta X_s \rangle
            \right\}
    \right\rVert       \notag \\
    & \qquad\quad +
    \left\lVert
        \sum_{ s \in D^{\varepsilon} \cap [0,t] } 
            D_a f(A_{s-},X_{s-})\Delta A_s
    \right\rVert
    + 
    \left\lVert
        \sum_{ s \in D^{\varepsilon} \cap [0,t] } 
            D_x^2 f(A_{s-},X_{s-})(\Delta X_s)^{\otimes 2} 
    \right\rVert      \notag \\
    & \qquad\quad +
    \varlimsup_{n \to \infty}
        \sup_{\substack{z,w \in K_1  \\ \lvert z-w \rvert \leq \alpha(\varepsilon,n)}}
            \rho_{K_3}(D_B^2 f(z) - D_B^2 f(w))
        V^2(X;\Pi)_t,   \notag \\
    & \qquad\quad +
    \varlimsup_{n \to \infty}
        \sup_{\substack{z,w \in K_1 \times K_2  \\ \lvert z-w \rvert \leq \alpha(\varepsilon,n) + \beta(\varepsilon,n)}}
        \rho_{K_1-K_1}(D_a f(z) - D_a f(w))
        V(A)_t.
\end{align}
Finally, by letting $\varepsilon \to 0$,
we see that the right-hand side of \eqref{5.9i} converges weakly to 
that of \eqref{5.9e}. This completes the proof for the weak case.

If $Q_B(X,X)$ is the strong $B$-quadratic variation,
then $I_5^{(n)}$ converges in the norm topology by the strong version of Lemma~\ref{2g}.
In this case, we obtain the norm convergence of the Riemannian sums
by replacing \eqref{5.9j} with a similar norm estimate.
\end{proof}

Combining Corollaries~\ref{2fab} and \ref{5.6d},
we obtain the integration by parts formula.
Note that the existence of the It\^{o}--F\"{o}llmer
integral $\int_{0}^{t} A_{s-} \otimes \mathrm{d}X_{s}$
follows from Corollary~\ref{2fab} and Proposition~\ref{5.9d}.

\begin{corollary}
Let $(\pi_n)$, $X$, and $A$ satisfy the same assumptions
as Corollary~\ref{2fab}. Then,
\begin{equation*}
    A_t \otimes X_t
    = 
    \int_{0}^{t} \mathrm{d}A_s \otimes X_{s-}
    + \int_{0}^{t} A_{s-} \otimes \mathrm{d}X_{s}
    + {^\alpha [A,X]}_t.
\end{equation*}
\end{corollary}

\appendix

\section{Differential calculus in Banach spaces}

In this section, we review some basic results about differential calculus in Banach spaces.

\begin{definition} \label{a1b}
Let $E$ and $F$ be Banach spaces and $U$ an open subset of $E$.
\begin{enumerate}
    \item A function $f \colon U \to F$ is G\^{a}teaux differentiable at $x \in U$
        if there exits an $L \in \mathcal{L}(E,F)$ such that
        \begin{equation*}
            \lim_{t \to 0, t \neq 0} \frac{1}{t} \{ f(x+th) - f(x) \} = Lh
        \end{equation*}
        for all $h \in E$.
        The function $f$ is G\^{a}teaux differentiable if it is 
        G\^{a}teaux differentiable at all points in $U$.
    \item A function $f \colon U \to F$ is Fr\'{e}chet differentiable at $x \in U$
        if there exits an $L \in \mathcal{L}(E,F)$ such that
        \begin{equation*}
            \lVert f(x+h) - f(x) - Lh \rVert = o(\lVert h \rVert).
        \end{equation*}
        The function $f$ is Fr\'{e}chet differentiable if it is 
        Fr\'{e}chet differentiable at all points in $U$.
\end{enumerate}
\end{definition}

The bounded operators $L$ in (i) and (ii) are called G\^{a}teaux and Fr\'{e}chet 
derivatives, respectively, and denoted by the same symbol $Df(x)$.
If $f$ is Fr\'{e}chet differentiable, then it is G\^{a}teaux differentiable
and both derivatives coincide.
Hence this notation is consistent.

Higher-order G\^{a}teaux and Fr\'{e}chet derivatives 
are defined inductively by the formula $D^{n+1} f = D D^n f$.
The $n$-th G\^{a}teaux derivative $D^{n}f(x)$ at $x \in U$ is an element of 
$\mathcal{L}(E^{\widehat{\otimes} n},F) \cong \mathcal{L}^{(n)}(E^n;F)$,
where $\mathcal{L}^{(n)}(E^n;F)$ denotes the set of bounded $n$-linear maps.
Moreover, if $f$ is $n$-times Fr\'{e}chet differentiable,
the multilinear map $D^n f(x) \colon E^n \to F$ is symmetric.

Before introducing Taylor's theorem,
we define a mild continuity condition
for a function defined on a subset of a Banach space.

\begin{definition} \label{a1bb}
Let $E$ be Banach spaces and $U$ an open convex subset of $E$.
A function $f \colon U \to T$ into a topological space is \emph{continuous along line segments} if, 
for every $x$ and $y$ in $U$, 
its restriction to the line segment $\{ \theta x + (1-\theta)y \mid \theta \in [0,1] \}$
is continuous.
\end{definition}

One can easily see that a G\^{a}teaux differentiable function $f \colon U \to F$
is continuous along line segments.

The following proposition gives a version of Taylor's formula 
that we use in this article.
It plays an essential role to prove the It\^{o} formula in Section 7.
Although that theorem seems classical, we give a proof for the reader's convenience.

\begin{proposition} \label{a1c}
Let $E$ and $F$ be Banach spaces and $U$ an open convex subset of $E$.
Assume that $f \colon U \to F$ is $n$-times G\^{a}teaux differentiable and
$D^n f$ is strongly continuous along line segments.
Then it admits Taylor's formula
\begin{equation*}
    f(x+u) - f(x)
    = 
    \sum_{1 \leq k \leq n-1} \frac{1}{k!} D^{k}f(x)u + 
    \int_0^1 \frac{(1-\theta)^{n-1}}{(n-1)!} D^n f(x+ \theta u) u^{\otimes n} \, \mathrm{d}\theta
\end{equation*}
for all $x \in U$ and $u \in E$ with $x+u \in U$.
\end{proposition}

\begin{proof}
The proof is by induction on $n$.

Suppose that $f$ is G\^{a}teaux differentiable and $Df$ is strongly continuous along line segments.
Then the map $\theta \mapsto f(x+ \theta u)$ has the continuous derivative
\begin{equation*}
    \frac{\mathrm{d}}{\mathrm{d} \theta} f(x+ \theta u) = Df(x+\theta u)u.
\end{equation*}
Hence, by the fundamental theorem of calculus, we get
\begin{equation*}
    f(x+u) - f(x) = \int_0^1 Df(x+ \theta u)u \, \mathrm{d}\theta.
\end{equation*}

Next, assume that the assertion holds for $n$.
Let $f \colon U \to F$ be an $n+1$-times G\^{a}teaux differentiable function
whose derivative $D^{n+1}f$ is strongly continuous along line segments.
By assumption, $f$ satisfies
\begin{equation} \label{a1e}
    f(x + u) - f(x) 
    = 
    \sum_{1 \leq k \leq n-1} \frac{1}{k!}D^{k}f(x)u^{\otimes k} 
    + \int_0^1 \frac{(1-\theta)^{n-1}}{(n-1)!} D^{n}f(x + \theta u)u^{\otimes n} \mathrm{d}\theta.
\end{equation}
Computation by the Leibniz rule shows
\begin{equation*}
    \frac{\mathrm{d}}{\mathrm{d}\theta} \left( \frac{(1-\theta)^{n}}{n!}D^{n}f(x + \theta u)u^{\otimes n} \right) 
    = -\frac{(1-\theta)^{n-1}}{(n-1)!} D^{n}f(x + \theta u)u^{\otimes n} +  \frac{(1-\theta)^{n}}{n!}D^{n+1}f(x + \theta u)u^{\otimes n+1}.
\end{equation*}
Therefore,
\begin{align} \label{a1d}
    - \frac{1}{n!}D^{n}f(x )u^{\otimes n} 
    = 
    - \int_0^1 \frac{(1-\theta)^{n-1}}{(n-1)!} D^{n}f(x + \theta u)u^{\otimes n} \mathrm{d}\theta 
    + \int_0^1 \frac{(1-\theta)^{n}}{n!}D^{n+1}f(x + \theta u)u^{\otimes n+1} \mathrm{d}\theta.
\end{align}
Combining \eqref{a1e} and \eqref{a1d}, we get the equation
\begin{equation*}
    f(x + u) - f(x) = \sum_{1 \leq k \leq n} \frac{1}{k!}D^{k}f(x)u^{\otimes k} + \int_0^1 \frac{(1-\theta)^{n}}{n!} D^{n+1}f(x + \theta u)u^{\otimes n+1} \mathrm{d}\theta.
\end{equation*}
This completes the proof.
\end{proof}

The following proposition is also used in the proof of the main theorem.

\begin{proposition} \label{a1f}
Let $E_1$, $E_2$, and $F$ be Banach spaces.
Assume that $f \colon E_1 \times E_2 \to F$ satisfies the following conditions:
\begin{enumerate}
    \item The map $x_1 \mapsto f(x_1,x_2)$ is G\^{a}teaux differentiable for all $x_2 \in E_2$ 
        and $D_{x_1} f \colon E_1 \times E_2 \to \mathcal{L}(E_1,F)$ is strongly continuous;
    \item The map $x_2 \mapsto f(x_1,x_2)$ is G\^{a}teaux differentiable for all $x_1 \in E_1$
        and $D_{x_2} f \colon E_1 \times E_2 \to \mathcal{L}(E_2,F)$ is strongly continuous.
\end{enumerate}
Then $f \colon E_1 \times E_2 \to F$ has the G\^{a}teaux derivative given by 
$Df = (D_{x_1}f,D_{x_2} f)$ and $Df \colon E_1 \times E_2 \to \mathcal{L}(E_1 \oplus E_2,F)$
is strongly continuous.
Moreover, $Df\vert_{K}$ induces a continuous map
$Df\vert_K \colon K \times (E_1 \times E_2) \to F$ for every compact set $K \subset E_1 \times E_2$.
\end{proposition}

\begin{proof}
Fix an $(x_1,x_2) \in E_1 \times E_2$
and take an arbitrary directional vector $(h_1,h_2) \in E_1 \times E_2$.
Applying Taylor's formula to the first variable, we get
\begin{equation*}
    f(x_1 + t h_1,x_2 + t h_2) - f(x_1,x_2)
    =
    \int_0^1 D_{x_1} f(x_1 + \theta t h_1,x_2 + t h_2) t h_1 \, \mathrm{d}\theta
    + f(x_1,x_2 + t h_2) - f(x_1,x_2)
\end{equation*}
for all $t \neq 0$.
Since $D_{x_1}f$ is strongly continuous, we see that
\begin{equation*}
    \lim_{t \to 0} D_{x_1} f(x_1 + \theta t h_1,x_2 + t h_2) h_1 = D_{x_1}f(x_1,x_2)h_1, \qquad
    \sup_{\theta, t \in \lrbrack 0,1 \rbrack} \lVert D_{x_1} f(x_1 + \theta t h_1,x_2 + t h_2) \rVert < \infty.
\end{equation*}
Note that the inequality above follows from the uniform boundedness principle.
Therefore, by the dominated convergence theorem,
\begin{equation*}
    \frac{1}{t} \int_0^1 D_{x_1} f(x_1 + \theta t h_1,x_2 + t h_2) t h_1 \, \mathrm{d}\theta
    = 
    \int_0^1 D_{x_1} f(x_1 + \theta t h_1,x_2 + t h_2) h_1 \, \mathrm{d}\theta
    \xrightarrow[t \to 0]{} f(x_1,x_2)h_1.
\end{equation*}
This shows 
\begin{equation*}
    \lim_{t \to 0, t\neq 0} \frac{1}{t} \{ f(x_1 + t h_1,x_2 + t h_2) - f(x_1,x_2) \}
    =
    D_{x_1}f(x_1,x_2)h_1 + D_{x_2}f(x_1,x_2)h_2,
\end{equation*}
which means that $(D_{x_1}f,D_{x_2}f)$ is the G\^{a}teaux derivative of $f$
at $(x_1,x_2)$.
The strong continuity of $Df = (D_{x_1}f,D_{x_2}f)$ is obvious.

The last claim follows from Lemma~\ref{a1g} below.
\end{proof}

\begin{lemma} \label{a1g}
Let $K$ be a compact topological space and let $E$ and $F$ be Banach spaces.
Assume that a map $\varphi \colon K \to \mathcal{L}(E,F)$ is strongly continuous.
Then, the map $\widetilde{\varphi} \colon K \times E \to F$
induced by $\varphi$ is continuous.
\end{lemma}

\begin{proof}
Fix $(t,x) \in K \times E$ and take 
an arbitrary net $(t_{\lambda},x_{\lambda})_{\lambda \in \Lambda}$
that converges to $(t,x)$ in the product topology.
Then 
\begin{equation*}
    \lVert \varphi(t)x - \varphi(t_{\lambda})x_\lambda \rVert 
    \leq 
    \lVert \varphi(t)x - \varphi(t_{\lambda})x \rVert 
    + \sup_{t \in K} \lVert \varphi(t) \rVert \lVert x - x_{\lambda} \rVert.
\end{equation*}
Because $\varphi$ is strongly continuous, the first term on the right-hand side converges to $0$.
The strong continuity and the uniform boundedness principle imply
that $\sup_{t \in K} \lVert \varphi(t) \rVert < \infty$.
Hence the second term also converges to $0$.
As a consequence, we find that $\widetilde{\varphi}(t_{\lambda},x_{\lambda}) \coloneqq \varphi(t_{\lambda})x_{\lambda}$
converges to $\widetilde{\varphi}(t,x) \coloneqq \varphi(t)x$.
\end{proof}

\section{Vector integration}

\subsection{A brief review of vector integration in Banach spaces}

In this subsection, we will review an integration theory for vector functions
and vector measures of finite variation.
Let $E$, $F$, and $G$ be three Banach spaces
and $B \colon F \times E \to G$ be a bounded bilinear map.
We aim to introduce integrals of the form $\int B(f,\mathrm{d}\mu)$,
where $f$ is an $F$-valued `nice' function and $\mu$ is an $E$-valued 
$\sigma$-additive measure defined on a $\delta$-ring of subsets of $\mathbb{R}_{\geq 0}$.
See Dinculeanu~\cite[Chapter 1 {\S} 2]{Dinculeanu_2000}
for details about the contents of this section.
We also refer to Diestel and Uhl~\cite{Diestel_Uhl_1977}
and Dinculeanu~\cite{Dinculeanu_1967},
which are classical references in the theory of vector measures and integration.

Let $\mathcal{I}$ be the semiring of subsets of $\mathbb{R}_{\geq 0}$
consisting of all bounded intervals of the form $\lrbrack a,b \rbrack$ and the singleton $\{ 0 \}$.
Moreover, let $\mathcal{D}$ be the $\delta$-ring generated by $\mathcal{I}$.
Note that both $\mathcal{I}$ and $\mathcal{D}$ generate the Borel $\sigma$-algebra.
The variation of a $\sigma$-additive vector measure $\mu \colon \mathcal{D} \to E$
on a subset  $A$ of $\mathbb{R}_{\geq 0}$ is defined by the formula
\begin{equation*}
    \lvert \mu \rvert(A)
    \coloneqq 
    \sup \left\{ 
        \sum_{\lambda \in \Lambda} \lVert \mu(A_\lambda) \rVert_E \,
        \middle\vert\, 
        \text{
            $(A_\lambda)_{\lambda \in \Lambda}$: a finite disjoint family of elements
            of $\mathcal{D}$
        }, \ 
        \bigcup_{\lambda \in \Lambda} A_\lambda \subset A
    \right\}.
\end{equation*}
We say that $\mu$ has \emph{finite variation} if $\lvert \mu \rvert(A) < \infty$
for all $A \in \mathcal{D}$.
The measure $\mu$ has \emph{bounded variation} whenever
$\lvert \mu \rvert(\mathbb{R}_{\geq 0}) < \infty$.
Since $\mu$ is assumed to be $\sigma$-additive, 
the variation measure $\lvert \mu \rvert \colon \mathcal{D} \to \mathbb{R}_{\geq 0}$
is also $\sigma$-additive. 
Then there is a unique $\sigma$-additive measure,
denoted by the same symbol $\lvert \mu \rvert$, 
on the $\sigma$-algebra $\mathcal{B}(\mathbb{R}_{\geq 0})$
that extends $\lvert \mu \rvert \colon \mathcal{D} \to \mathbb{R}_{\geq 0}$.

An $F$-valued $\mathcal{D}$-simple function is a function $f \colon \mathbb{R}_{\geq 0} \to F$ of the form 
\begin{equation*}
    f = \sum_{\lambda \in \Lambda} 1_{A_\lambda} a_\lambda ,
\end{equation*}
where $\Lambda$ is a finite set,
$(a_\lambda)_{\lambda \in \Lambda} \in E^\Lambda$,
and $(A_\lambda)_{\lambda \in \Lambda} \in \mathcal{D}^\Lambda$
is a disjoint family.
We can immediately define the integral of $f$ by
\begin{equation*}
    \int_{\mathbb{R}_{\geq 0}} B(f,\mathrm{d}\mu)
    =
    \sum_{\lambda \in \Lambda} B(a_\lambda,\mu(A_\lambda)).
\end{equation*}
Let $S_F(\mathcal{D})$ denote the set of all 
$F$-valued $\mathcal{D}$-simple functions.
Then $S_F(\mathcal{D})$ is dense in the Lebesgue--Bochner space 
\begin{equation*}
    L^1(\lvert \mu \rvert;F) 
    \coloneqq
    \left\{ 
        f \colon \mathbb{R}_{\geq 0} \to F  \,
        \middle\vert \,
        \text{$f$ is strongly $\mathcal{B}(\mathbb{R}_{\geq 0})$-measurable, }
        \int_{\mathbb{R}_{\geq 0}} \lVert f(s) \rVert_F \, \lvert \mu \rvert (\mathrm{d}s) < \infty 
    \right\}.
\end{equation*}
Since the integration map $f \mapsto \int_{\mathbb{R}_{\geq 0}} B(f,\mathrm{d}\mu)$
is bounded and linear, it can be uniquely extended to a bounded linear operator
$T_{\mu} \colon L^1(\lvert \mu \rvert;F) \to G$.
If a strongly measurable function $f \colon \mathbb{R}_{\geq 0} \to F$
and an $A \in \mathcal{B}(\mathbb{R}_{\geq 0})$ satisfy $f1_A \in L^1(\lvert \mu \rvert;F)$,
then we define the integral of $f$ on $A$ as
\begin{equation*}
    \int_{A} B(f(s),\mu(\mathrm{d}s))
    =
    \int_{A} B(f,\mathrm{d}\mu)
    \coloneqq
    T_{\mu}(f1_A).
\end{equation*}
In particular, every $f \in L^{1}_{\mathrm{loc}}(\lvert \mu \rvert;F)$
can be integrated on every bounded interval. 
We also use the notation 
\begin{equation*}
    \int_{A} B(f, \mathrm{d}\mu)
    =
    \int_a^b B(f, \mathrm{d}\mu)
\end{equation*}
if $A$ is a bounded interval of the form $A = \lrbrack a,b \rbrack$.
By a direct calculation, one can derive the inequality
\begin{equation*}
    \left\lVert \int_{A} B(f,\mathrm{d}\mu) \right\rVert_G
    \leq 
    \lVert B \rVert \int_{A} \lVert f(s) \rVert_E \lvert \mu \rvert(\mathrm{d}s)
\end{equation*}
for all such $A$ and $f$.
This estimate guarantees that 
the dominated convergence theorem remains valid in this situation.

Finally, we introduce a decomposition of
a vector measure on $\mathbb{R}_{\geq 0}$
into atomic and nonatomic parts.
As before, let $\mu \colon \mathcal{D} \to E$ be a $\sigma$-additive measure
of finite variation. Set
\begin{equation*}
    D = \{ t \in \mathbb{R}_{\geq 0} \mid \lvert \mu \rvert (\{ t \}) > 0 \}.
\end{equation*}
Because $\mu$ has finite variation, we see that $D$ is countable.
Besides, since each singleton $\{ t \}$ belongs to $\mathcal{D}$,
we have the equality $\lVert \mu(\{ t \}) \rVert = \lvert \mu \rvert (\{ t \})$
for all $t \in \mathbb{R}_{\geq 0}$.
Now let
\begin{equation*}
    \mu^{\mathrm{d}} = \sum_{s \in D} \delta_s \mu(\{ s \}),
\end{equation*}
where $\delta_s$ denotes the Dirac measure at $s \in \mathbb{R}_{\geq 0}$.
Then $\mu^{\mathrm{d}}$ is a $\sigma$-additive vector measure of finite variation 
that satisfies 
\begin{equation*}
    \lvert \mu^{\mathrm{d}} \rvert(A) = \sum_{s \in D} \lVert \mu(\{ s \}) \rVert_E \, \delta_s(A) < \infty
\end{equation*}
for each $A \in \mathcal{D}$.
If we define
$\mu^{\mathrm{c}} = \mu - \mu^{\mathrm{d}}$,
then $\mu^{\mathrm{c}}$ and $\mu^{\mathrm{d}}$ give a mutually singular
decomposition of $\mu$ satisfying
$\lvert \mu \rvert = \lvert \mu^{\mathrm{c}} \rvert + \lvert \mu^{\mathrm{d}} \rvert$.
This decomposition of a measure gives a decomposition of an integral as
\begin{equation*}
    \int_{\mathbb{R}_{\geq 0}} B(f,\mathrm{d}\mu)
    = 
    \int_{\mathbb{R}_{\geq 0}} B(f,\mathrm{d}\mu^{\mathrm{c}})
    + \int_{\mathbb{R}_{\geq 0}} B(f,\mathrm{d}\mu^{\mathrm{d}})
\end{equation*}
for every $f \in L^1(\lvert \mu \rvert;F)$.

\subsection{Extension of vector integration}

In Section B.1, integrands are assumed to be strongly measurable
with respect to the norm topology of $F$;
i.e. integrands can necessarily be approximated pointwise by simple functions 
in the norm topology.
A function $f \colon \mathbb{R}_{\geq 0} \to F$
that is c\`{a}dl\`{a}g in a weaker topology may not satisfy this condition.
We need, however, to integrate such functions for our It\^{o}--F\"{o}llmer formula
(Theorem~\ref{2e}).
For this purpose, we extend the vector integration introduced in 
a previous section to a suitable setting.

As in Appendix B.1, 
let $E$ and $G$ be Banach spaces and
let $\mu \colon \mathcal{D} \to E$
be a $\sigma$-additive vector measure of finite variation.

\begin{lemma} \label{a3b}
The range $\mu(\mathcal{I} \cap [0,T])$ is relatively compact in $E$ for all $T > 0$. 
\end{lemma}

\begin{proof}
First, define a function $F \colon \mathbb{R}_{\geq 0} \to E$ by the formula
$F(t) = \mu([0,t])$. $\sigma$-additivity of $\mu$ implies that 
$F$ is c\`{a}dl\`{a}g in $E$.
For a given $T > 0$, we can find a compact set $K$ including the image $F([0,T])$.
If $I = \lrbrack a,b \rbrack \subset [0,T]$, then 
\begin{equation*}
    \mu(\lrbrack a,b \rbrack) = F(b) - F(a) \in K - K.
\end{equation*}
On the other hand, if $I = \{ 0 \} \subset [0,T]$, 
\begin{equation*}
    \mu(\{ 0 \}) = F(0) \in K.
\end{equation*}
Thus $\mu(\mathcal{I}) \subset K \cup (K - K)$.
Because $K$ and $K-K$ are both compact, $\mu(\mathcal{I})$ is relatively compact.
\end{proof}

\begin{lemma} \label{a3bb}
Let $F = \mathcal{L}(E,G)$ and let 
$S_F(\mathcal{I})$ be the set of all $F$-valued $\mathcal{I}$-simple functions.
Suppose that $T > 0$ and $K$ is a compact set 
such that $\mu(\mathcal{I} \cap [0,T]) \subset K$.
Then 
\begin{equation*}
    \left\lVert \int_{[0,T]} f \mathrm{d}\mu \right\rVert_{G}
    \leq 
    \sup_{s \in [0,T]} \rho_K(f(s)) \lvert \mu \rvert([0,T]). 
\end{equation*}
\end{lemma}

Here, recall that a seminorm $\rho_K$ is defined by \eqref{2.1b} in Section 2.1.

\begin{proof}
Let $f = \sum_{I \in \Lambda} 1_I a_I$ be an $\mathcal{I}$-simple function
with $\Lambda$ being disjoint. Then, by the definition of the integral and the variation of $\mu$,
we see that 
\begin{equation*}
    \left\lVert \int_{[0,T]} f \mathrm{d}\mu \right\rVert_{G}
    \leq 
    \sum_{I \in \Lambda} \rho_K(a_I) \lVert \mu(I \cap [0,T]) \rVert_E 
    \leq 
    \sup_{s \in [0,T]} \rho_K(f(s)) \lvert \mu \rvert ([0,T]).
\end{equation*}
This is the desired inequality.
\end{proof}

In what follows, let $D^-(\mathbb{R}_{\geq 0}, \mathcal{L}_{\mathrm{c}}(E,G))$ 
stands for the locally convex Hausdorff topological vector space 
consisting of all c\`{a}gl\`{a}d functions 
endowed with the topology of uniform convergence on compacta.

\begin{lemma} \label{a3bc}
The space $S_F(\mathcal{I})$ is dense in 
$D^-(\mathbb{R}_{\geq 0}, \mathcal{L}_{\mathrm{c}}(E,G))$.
\end{lemma}

\begin{proof}
First, let $\mathcal{K}$ the family of all compact subsets of $E$
and define a family of seminorms $(\rho_{T,K}; T > 0,K \in \mathcal{K})$ by the formula
\begin{equation*}
    \rho_{T,K}(f) = \sup_{s \in [0,T]} \rho_K(f(s)).
\end{equation*}
Then the topology of $D^-(\mathbb{R}_{\geq 0}, \mathcal{L}_{\mathrm{c}}(E,G))$
is induced by $(\rho_{T,K})$.

Set $\Lambda = \mathbb{R}_{> 0} \times \mathcal{K} \times \mathbb{R}_{> 0}$
and define an order on $\Lambda$ by letting  
$(T_1,K_1,\varepsilon_1) \leq (T_2,K_2,\varepsilon_2)$
whenever $T_1 \leq T_2$, $K_1 \subset K_2$, and $\varepsilon_1 \geq \varepsilon_2$.
This order makes $\Lambda$ a directed set 
because $(T_1 \vee T_2, K_1 \cup K_2, \varepsilon_1 \wedge \varepsilon_2) \in \Lambda$.
For each $\lambda = (T,K,\varepsilon)$, we can find a $\mathcal{I}$-simple function 
$s_{\lambda}$ satisfying $\rho_{T,K} (s_{\lambda} - f) < \varepsilon$
as in the proof of Lemma~\ref{5.7g}.
The net $(s_{\lambda})_{\lambda \in \Lambda}$ converges to $f$ in the topology of 
$D^-(\mathbb{R}_{\geq 0}, \mathcal{L}_{\mathrm{c}}(E,G))$.
Indeed, take $\varepsilon_0 > 0$ and $(T_0,K_0)$ arbitrarily.
If $\lambda = (T,K,\varepsilon) \geq (T_0,K_0,\varepsilon_0)$, then 
\begin{equation*}
    \rho_{T_0,K_0}(s_\lambda - f) 
    \leq 
    \rho_{T,K}(s_\lambda -f)
    < \varepsilon
    \leq \varepsilon_0.
\end{equation*}
Hence $(s_{\lambda})_{\lambda \in \Lambda}$ converges to $f$ in $\rho_{T_0,K_0}$.
Since the choice of $(T_0,K_0)$ is arbitrary, we can conclude that 
$S_F(\mathcal{I})$ is dense in 
$D^-(\mathbb{R}_{\geq 0}, \mathcal{L}_{\mathrm{c}}(E,G))$.
\end{proof}

\begin{theorem} \label{a3c}
Let $T > 0$. The integration map
\begin{equation*}
    S_F(\mathcal{I}) \ni f \longmapsto \int_{[0,T]} f \mathrm{d}\mu \in G
\end{equation*}
can be uniquely extended to a continuous linear map
on $D^-(\mathbb{R}_{\geq 0}, \mathcal{L}_{\mathrm{c}}(E,G))$.
\end{theorem}

\begin{proof}
Define the integral map $I_T \colon \mathcal{S}_F(\mathcal{I}) \to G$
by $I_T(f) = \int_{[0,T]}f \mathrm{d}\mu$.
Because $\mathcal{S}_F(\mathcal{I})$ is dense in $D^-(\mathbb{R}_{\geq 0}, \mathcal{L}_{\mathrm{c}}(E,G))$
by Lemma~\ref{a3bc},
it suffices to show that $I_T$ is linear and continuous.
Refer to Schaefer and Wolff~\cite[Chapter III Section 1]{Schaefer_Wolff_1999}
for the extension of a continuous linear map between topological vector spaces.

The linearity of $I_T$ is obvious.
The continuity follows from Lemma~\ref{a3bb} and 
a standard continuity criterion for linear maps in
topological vector spaces.
See, e.g. Schaefer and Wolff~\cite[Chapter III Section 1 (1.1)]{Schaefer_Wolff_1999}.
\end{proof}

We will show a variant of the dominated convergence theorem for our vector integrals.

\begin{lemma} \label{a3d}
For each $f \in D^-(\mathbb{R}_{\geq 0}, \mathcal{L}_{\mathrm{c}}(E,G))$
and compact subset $K$ of $E$,
the function $t \mapsto \rho_K(f(t))$ is Borel measurable.
\end{lemma}

\begin{proof}
If $f \in D^-(\mathbb{R}_{\geq 0}, \mathcal{L}_{\mathrm{c}}(E,G))$, then 
the map $t \mapsto \rho_K(f(t))$ is c\`{a}gl\`{a}d.
Hence it is Borel measurable.
\end{proof}

\begin{lemma} \label{a3e}
Let $f = \sum_{I \in \Lambda} 1_{I} a_{I}$ 
be an $\mathcal{L}_\mathrm{c}(E,G)$-valued $\mathcal{I}$-simple function
such that $\Lambda$ is disjoint.
Then,
\begin{equation*}
    \rho_K(f(t)) = \sum_{I \in \Lambda} \rho_K(a_I) 1_{I}(t), \qquad t \geq 0.
\end{equation*} 
\end{lemma}

\begin{proposition} \label{a3f}
Let $T > 0$ and let $K$ be a compact set satisfying $\mu(\mathcal{I} \cap [0,T]) \subset K$.
Then 
\begin{equation} \label{a3fb}
    \left\lVert \int_{[0,T]} f \mathrm{d}\mu \right\rVert_G
    \leq
    \int_{[0,T]} \rho_K(f(s)) \lvert \mu \rvert(\mathrm{d}s)
\end{equation}
for all $f \in D^-(\mathbb{R}_{\geq 0}, \mathcal{L}_{\mathrm{c}}(E,G))$.
\end{proposition}

\begin{proof}
If a simple function $f$ has a disjoint representation
$f = \sum_{I \in \Lambda} 1_I a_I$,
then
\begin{align*}
    \left\lVert \int_{[0,T]} f \mathrm{d}\mu \right\rVert_{G}
    & \leq 
    \sum_{I \in \Lambda} \rho_K(a_I) \lVert \mu(I \cap [0,T]) \rVert_E   \\
    & \leq 
    \sum_{I \in \Lambda} \rho_K(a_I) \lvert \mu \rvert (I \cap [0,T])   \\
    & = 
    \int_{[0,T]} \rho_K(f(s)) \lvert \mu \rvert(\mathrm{d}s).
\end{align*}
Note that the last equality follows from Lemma~\ref{a3e}.
Therefore \eqref{a3fb} holds on $S_F(\mathcal{I})$.

The general case is proved by approximation.
\end{proof}

\begin{proposition}[Dominated convergence theorem] \label{a3g}
Let $T > 0$ and let $K$ be a compact subset of $E$ $\mu(\mathcal{I} \cap [0,T]) \subset K$.
Suppose that a sequence $(f_n)$ and an element $f$ 
in $D^-(\mathbb{R}_{\geq 0}, \mathcal{L}_{\mathrm{c}}(E,G))$
satisfy the following conditions:
\begin{enumerate}
    \item The sequence $(f_n)$ converges to $f$ pointwise on $[0,T]$ with respect to $\rho_K$;
    \item There is a $g \in L^1([0,T],\lvert \mu \rvert)$
        such that $\rho_K(f_n(t)) \leq g(t)$ almost everywhere on $[0,T]$. 
\end{enumerate}
Then,
\begin{equation*}
    \left\lVert \int_{[0,T]} f_n \mathrm{d}\mu - \int_{[0,T]} f \mathrm{d}\mu \right\rVert_G
    \xrightarrow[n \to \infty]{} 0.
\end{equation*}
\end{proposition}

\begin{proof}
Proposition~\ref{a3f} implies the inequality
\begin{equation*}
    \left\lVert \int_{[0,T]} f_n \mathrm{d}\mu - \int_{[0,T]} f \mathrm{d}\mu \right\rVert_G
    \leq 
    \int_{[0,T]} \rho_K(f_n(s) - f(s)) \lvert \mu \rvert (\mathrm{d}s).
\end{equation*}
By applying the dominated convergence theorem to 
the integral on the right-hand side,
we obtain the desired convergence.
\end{proof}

\begin{proposition} \label{a3h}
If $f \in D^-(\mathbb{R}_{\geq 0}, \mathcal{L}_{\mathrm{c}}(E,G))$, then
\begin{equation*}
    \left\langle z^*, \int_{\lrbrack 0,t \rbrack} f(s) \mu(\mathrm{d}s) \right\rangle
    =
    \int_{\lrbrack 0,t \rbrack} \langle z^*f(s), \mu(\mathrm{d}s) \rangle
\end{equation*}
for all $z^* \in G^*$ and $t \in \mathbb{R}_{\geq 0}$.
\end{proposition}

\begin{proof}
First, assume that $f$ has the form
$f = \sum_{I \in \Lambda} 1_{I} a_{I}$ with $\Lambda \subset \mathcal{I}$ disjoint.
By a direct calculation, we see that
\begin{equation*}
    \left\langle z^*, \int_{\lrbrack 0,t \rbrack} f(s) \mu(\mathrm{d}s) \right\rangle
    =
    \sum_{I \in \Lambda} z^* a_I \mu(I \cap \lrbrack 0,t])
    =
    \int_{\lrbrack 0,t \rbrack} \left\langle z^* f(s), \mu(\mathrm{d}s) \right\rangle.
\end{equation*}
Hence the formula holds for simple functions.
It can be extended to general integrands by the density argument. 
\end{proof}

\paragraph{Acknowledgements}
The author thanks 
Professor Jun Sekine, his doctoral supervisor, and Professor Masaaki Fukasawa
for their helpful comments and encouragement.
The author is also grateful to the anonymous referees 
for their various useful comments on an earlier version of this paper.
This work was partially supported by
JSPS KAKENHI Grant Number JP19H00643 (Masanori Hino).

\printbibliography[heading=bibintoc]

Graduate School of Science,
Kyoto University, Kyoto, Japan

\emph{Current address}:
Department of Creative Engineering, National Institute of Technology, Tsuruoka College, Japan

\emph{Email address}: \texttt{hirai@tsuruoka-nct.ac.jp}
\end{document}